\documentclass[a4paper,11pt]{amsart}
\usepackage[latin1]{inputenc}
\usepackage{amsmath}
\usepackage{amsfonts}
\usepackage{amsthm}

\numberwithin{equation}{section} 

\renewcommand{\t}{{\tau}}
\newcommand{\ep}{{\epsilon}}
\renewcommand{\l}{{\lambda}}
\newcommand{\G}{{\mathbb G}}
\newcommand{\cG}{{\mathcal G}}
\newcommand{\cS}{{\mathcal S}}

\newcommand{\ts}{{\tau_{\Sigma}}}
\newcommand{\Sxr}{{\Sigma_{x,r}}}
\newcommand{\f}{{\varphi}}
\newcommand{\g}{{\gamma}}
\newcommand{\F}{{\Phi}}

\newcommand{\R}{\mathbb R}
\newcommand{\haus}{{\mathcal H}}
\newcommand{\Shaus}{{\mathcal S}}
\newcommand{\lan}{\langle}
\newcommand{\ran}{\rangle}
\newcommand{\res}{\mathop{\hbox{\vrule height 7pt width .5pt depth 0pt
\vrule height .5pt width 6pt depth 0pt}}\nolimits}
\newcommand{\hsd}[1]{{\mathcal S}_d^{1}}
\newcommand{\hs}[2]{S^+_{\H}({1},{2})}
\newcommand{\N}{\mathbb N}
\newtheorem{teo}{Theorem}[section]
\newtheorem{lem}[teo]{Lemma}

\newtheorem{cor}[teo]{Corollary}
\newtheorem{defi}[teo]{Definition}
\newtheorem{rem}[teo]{Remark}
\newtheorem{ex}[teo]{Example}

\begin{document}

%%%%%%%%%%%%%%%%%%%%%%%%%%%%%%%%%%%%%%%%%%%%%%%%%%%%%%%%%%%%%%%%%%%%%%%%

\title{An intrinsic measure for submanifolds in stratified groups}
\author{Valentino Magnani}
\address{Valentino Magnani: Dipartimento di Matematica
\\ Largo Bruno Pontecorvo 5 \\ 56127 Pisa - Italy}
\email{magnani@dm.unipi.it}
\author{Davide Vittone}
\address{Davide Vittone: Scuola Normale Superiore
\\ Piazza dei Cavalieri 7\\ 56126 Pisa - Italy} 
\email{d.vittone@sns.it}

%%%%%%%%%%%%%%%%%%%%%%%%%%%%%%%%%%%%%%%%%%%%%%%%%%%%%%%%%%%%%%%%%%%%%%%%%%

%\date{\today}
\maketitle

%%%%%%%%%%%%%%%%%%%%%%%%%%%%%%%%%%%%%%%%%%%%%%%%%%%%%%%%%%%%%%%%%%%%%%%%%5

\begin{quote}
{\textsc{Abstract.}}
For each submanifold of a stratified group,
we find a number and a measure only depending on
its tangent bundle, the grading and the fixed 
Riemannian metric. In two step stratified groups,
we show that such number and measure coincide
with the Hausdorff dimension and with
the spherical Hausdorff measure of the submanifold
with respect to the Carnot-Carath\'eodory distance,
respectively.
Our main technical tool is an intrinsic blow-up at
points of maximum degree.
We also show that the intrinsic tangent cone
to the submanifold at these points is always a subgroup. 
Finally, by direct computations in the Engel group,
we show how our results can be extended to higher step
stratified groups, provided the submanifold is sufficiently regular.

\end{quote}

\tableofcontents

\pagebreak

\section{Introduction}

In this paper we study how a submanifold inherits its sub-Riemannian geometry
from a stratified group equipped with its Carnot-Carath\'eodory distance.
Our aim is finding the sub-Riemannian measure
``naturally'' associated with a submanifold.

This measure for hypersurfaces
is exactly the $\mathbb G$-perimeter,
which is widely acknowledged as the appropriate measure
in connection with intrinsic regular hypersurfaces,
trace theorems, isoperimetric inequalities,
the Dirichlet problem for sub-Laplacians, minimal surfaces, 
and more. Here we address the reader 
to some relevant papers
\cite{Amb1}, \cite{ASCV}, \cite{Bal}, \cite{BSCV}, \cite{CDG}, \cite{CG},
\cite{CheHwa}, \cite{CHMY}, \cite{CitMan},
\cite{CitSar}, \cite{DGN}, \cite{DGN1}, \cite{DGN3},
\cite{DGN4}, \cite{FGW}, \cite{GN}, \cite{GarPau},
\cite{HlaPau}, \cite{KirSer}, \cite{FSSC4},
\cite{FSSC5}, \cite{LeoMas}, \cite{LeoRig},
\cite{Mag2}, \cite{Mag5}, \cite{Mag8}, \cite{Mat05}, \cite{MonMor1},
\cite{Pan01}, \cite{Pau3}, \cite{RitRos}, \cite{RitRos1}
and the reference therein.

Our question is: what does replace the $\mathbb G$-perimeter
in arbitrary submanifolds?
Clearly, once the Hausdorff dimension of the submanifold is known,
the corresponding spherical Hausdorff measure
should be the natural candidate. 
However this measure is not manageable, since it
cannot be used in minimization problems,
due to the lack of lower semincontinuity with respect to the
Hausdorff convergence of sets.
It is then convenient to find an equivalent measure,
that can be represented as the supremum among a
suitable family of linear functionals,
in analogy with the classical theory of currents.

In the recent works \cite{FSSC6}, \cite{Mag7}, higher codimensional
submanifolds in the Heisenberg group have been considered along with
their associated measure. For regular submanifolds, in \cite{Mag8}
this study is developped for the class of non-horizontal
submanifolds, as we will explain below.
Here we emphasize examples of H\"older submanifolds where
the Hausdorff measure with respect to the
Carnot-Carath\'eodory distance is finite, but the Riemannian 
measure is not, \cite{KirSer}.
Never\-theless, in \cite{FSSC6} the authors consider
intrinsic currents in the Heisenberg groups that include 
the previously mentioned ``singular'' submanifolds.

Our intrinsic measure has been found in \cite{Mag8} 
for the class of {\em non-horizontal submanifolds}, characterized
by having positive $(Q$$-$$k)$-dimensional spherical
Hausdorff measure $\cS^{Q-k}$.
Here $Q$ denotes the homogeneous dimension of the group and 
$k$ is the codimension of the submanifold.
Results of \cite{Mag8} rely on the following facts.
The Riemannian surface measure admits an intrinsic
blow-up at non-horizontal points and the other points
form an $\cS^{Q-k}$-negligible subset, according to \cite{Mag5}.
As a continuation of \cite{Mag8}, we wish to find the intrinsic
measure associated with 
the remaining manifolds, namely, {\em horizontal submanifolds}.
To do this, we have to find out the privileged
subset of points of a horizontal submanifold
where the blow-up holds.

Recall that at a horizontal point $x$ of a $C^1$ smooth
submanifold $\Sigma$ contained in a stratified group $\G$, 
the horizontal subspace $H_x\G$ and the tangent space $T_x\Sigma$
do not span all of $T_x\G$. 
We say that a submanifold is {\em horizontal} if it is formed
by horizontal points and {\em non-horizontal} otherwise.
Recall that horizontal points of hypersurfaces
coincide with the well known characteristic points, 
that play an important role in the study of hypersurfaces
in stratified groups, \cite{Bal}, \cite{CG}, \cite{DGN},
\cite{DGN1}, \cite{Der}, \cite{FraWhe}, \cite{FSSC4}, \cite{FSSC5},
\cite{Mag5}, \cite{MonMor}, \cite{MonMor2}.

Any smooth hypersurface is clearly non-horizontal,
due to the non-integrability of the horizontal distribution.
This is clearly not true in higher codimension,
where different situations can occur.
For instance, in the Heisenberg group $\mathbb H^n$ it is easy to check
that horizontal submanifolds exactly coincide with the special class of Legendrian submanifolds and it is easy to construct
non-horizontal submanifolds of any dimension.
On the other hand, there exist stratified groups where
all submanifolds of fixed topological dimension are horizontal,
see Example~\ref{exanothor}.

We first notice that horizontal points
may induce different behaviors of the submanifold when
it is dilated around these points.
We will show that this behavior depends on the {\em degree} $d_\Sigma(x)$
of the point $x$ in the submanifold $\Sigma$, see \eqref{dsx}
for precise definition.
This notion allows us to distinguish the different natures of horizontal 
points.
Roughly speaking, it represents a sort of ``pointwise Hausdorff dimension''.
Notice that our notion of degree for hypersurfaces satisfies the formula
$d_\Sigma(x)=Q-\mbox{{\em type}}(x)$, where the {\em type} of
a point in a hypersurface has been introduced in \cite{CG}.

The notion of degree permits us to characterize a
horizontal point $x\in\Sigma$, requiring that $d_\Sigma(x)<Q$$-$$k$.
At these points the blow-up of the submanifold, if it exists, it is not necessarily  a subgroup of $\G$, see Remark~\ref{nrpsubgroup}.
However, defining
$$
d(\Sigma)=\max_{x\in\Sigma}d_\Sigma(x)
$$
as the {\em degree} of $\Sigma$, we will
show that the blow-up always exists at points with maximum degree
$d_\Sigma(x)=d(\Sigma)$ and it is a subgroup of $\G$.
We have the following
\begin{teo}\label{blowup}
Let $\Sigma$ be a $C^{1,1}$ smooth submanifold of $\G$ and let
$x\in\Sigma$ be a point of maximum degree. Then for every $R>0$ we have  
\begin{equation}\label{blowupintro}
\delta_{1/r}(x^{-1}\Sigma)\cap D_R\rightarrow
\Pi_\Sigma(x)\cap D_R\quad\text{as }\quad r\rightarrow
 0^+
\end{equation}
with respect to the Hausdorff distance and $\Pi_\Sigma(x)$
is a subgroup of $\mathbb G$.
\end{teo}
Recall that $\delta_r$ are the intrinsic dilations of the group
and that $D_R$ is the closed ball of center the identity of the group
and radius $R$ with respect to the fixed homogeneous distance, see 
Section~\ref{prel} for details. 
The limit set $\Pi_\Sigma(x)$ corresponds to the one
introduced in Definition~\ref{defpix}.
In particular, Theorem~\ref{blowup} shows that the
{\em intrinsic tangent cone} to $\Sigma$ at $x$ exists,
according to Definition~3.4 of \cite{FSSC6},
and that it is exactly equal to $\Pi_\Sigma(x)$.

The geometrical interpretation of our approach
consists in foliating a neighbourhood of the point $x$
in $\Sigma$ with a family of curves which are homogeneous
with respect to dilations, up to infinitesimal terms of higher order.
In mathematical terms, we are able to represent $\Sigma$ in a neighbourhood
of $x$ as the union of curves $t\rightarrow\gamma(t,\lambda)$ in $\Sigma$
satisfying the Cauchy problem (\ref{gammalambda}).
These curves have
 the property
\begin{eqnarray}\label{dilparam}
	\gamma(t,\lambda)=\delta_t\left(G(\lambda)+O(t)\right)\,,
\end{eqnarray}
where $\lambda$ varies in a fixed compact set of $\R^p$
and the diffeomorphism $G$ defined in \eqref{defg}
parametrizes $\Pi_\Sigma(x)$ by $\mathbb R^p$
with respect to the graded coordinates, 
see Remark~\ref{remG}.
Our key tool is Lemma~\ref{curve} that shows
the crucial representation \eqref{dilparam} of the curves
parametrizing the submanifold.
The proof of this lemma is in turn due to the technical
Lemma~\ref{campisottoalg}, which is available since $\Pi_\Sigma(x)$
is a subgroup of $\G$.
From \eqref{blowupintro} we obtain the following
\begin{teo}\label{densita}
Let $\Sigma$ be a $C^{1,1}$ smooth $p$-dimensional submanifold of degree
$d=d(\Sigma)$ 
and let $x\in\Sigma$ be of the same degree. Then we have
\begin{equation}\label{eqdensita}
\lim_{r\downarrow 0}
\frac{\mbox{\rm vol}_{\tilde g}(\Sigma\cap B_{x,r})}{r^d} =
\frac{\theta(\tau_{\Sigma}^d(x))}{|\tau_{\Sigma}^d(x)|}\,.
\end{equation}
\end{teo}
Recall that $\theta(\tau_{\Sigma}^d(x))$ is the metric factor defined
in \eqref{metricfactor}, which also depends on the homogeneous distance
we are using to construct $\cS^d$. The simple $p$-vector
$\tau^d_\Sigma(x)$ is the part of unit tangent vector $\tau_\Sigma(x)$
having degree $d$, see \eqref{deftsm}. It generalizes the vertical
tangent vector introduced in \cite{Mag7} and \cite{Mag8}.
By \eqref{eqdensita} and standard theorems on differentiation
of measures, \cite{Fed}, we immediately deduce the following
\begin{equation}\label{integrfor}
\int_{\Sigma} \theta(\tau_{\Sigma}^d(x))\, d\Shaus^d_\rho(x)
= \int_{\Sigma} |\tau_{\Sigma}^d(x)|\,d\mbox{\rm vol}_{\tilde g}(x).
\end{equation}
whenever 
\begin{eqnarray}\label{dneglig}
\cS^d(\Sigma\setminus \Sigma_d)=0\,,
\end{eqnarray}
where $\Sigma_d$ is the open subset of points of maximum degree $d$.
In fact, it is not difficult to check that
$\tau^d_{\Sigma}$ vanishes on $\Sigma\setminus\Sigma_d$.
Formula \eqref{integrfor} shows that
$\cS^d$ is positive and finite on open bounded sets of the submanifold
and yields 
the ``natural'' sub-Riemannian measure on $\Sigma$:
\begin{eqnarray}\label{natural}
\mu=|\tau_{\Sigma}^d(x)|\;\mbox{\rm vol}_{\tilde g}\res \Sigma\,.
\end{eqnarray}
We stress that the measure defined in \eqref{natural}
does not depend on the metric $\tilde g$. In fact, parametrizing
a piece of $\Sigma$ by a mapping $\Phi:U\longrightarrow\G$, we have
$$
\mu\big(\Phi(U)\big)=\int_U
\left|\left(\partial_{x_1}\Phi\wedge
\partial_{x_2}\Phi\wedge\cdots\wedge\partial_{x_p}\Phi\right)_d\right|\,dx
$$   
where the projection $(\cdot)_d$ is defined in \eqref{taur}
and $|\cdot|$ is the norm induced by the fixed left invariant metric $g$.
This integral formula can be seen as an area-type formula
where the jacobian is projected on vectors of fixed degree.
From \eqref{integrfor} and the fact that $\theta(\cdot)$ is uniformly
bounded from below and from above, one easily deduces that $\mu$
is the ``natural'' replacement of $\cS^d$ and that it is a convenient
choice to solve sub-Riemannian filling problems, \cite{Gr1}.

In the case $d(\Sigma)=Q-k$, the limit \eqref{eqdensita} and
formula \eqref{integrfor} fit with (17) and (39) of \cite{Mag8} 
where $C^1$ smoothness of $\Sigma$ suffices, as we explain
in Remark~\ref{coin1cod}.
The validity of $d$-negligibility \eqref{dneglig} has been proved in
\cite{Mag5}, when $d(\Sigma)=Q$$-$$k$.
One can check that in two step stratified groups the formula
$2p-\dim(T_x\Sigma\cap H_x\G)=d_\Sigma(x)$ holds, then
estimates (3) of \cite{Mag8} immediately show that $d$-negligibility
holds in two step groups for submanifolds of arbitrary degree.
As a consequence,
$d(\Sigma)$ is the Hausdorff dimension of $\Sigma$.
However $d$-negligibility remains an interesting open question in
stratified groups of step higher than two, when $d(\Sigma)<Q$$-$$k$.

In the last part of this work
we study some examples of 2-dimensional submanifolds
of different degrees in the Engel group.
Despite $d$-negligibility  is an open question in groups of step
higher than two, our formula \eqref{integrfor} shows
the validity of \eqref{dneglig} for these examples.
This fact suggests that $d$-negligibility
should hold in any stratified group for submanifolds of arbitrary
degree, possibly requiring higher regularity.

\section{Preliminaries}\label{prel}
A stratified group $\G$ with topological dimension $q$ 
is a simply connected nilpotent Lie group with Lie algebra $\cG$
having the grading
\begin{equation}\label{stratificazione}
\cG=V_1\oplus\dots\oplus V_\iota\,,
\end{equation}
that satisfies the conditions $V_{i+1}=[V_1,V_i]$ for every $i\geq 1$ and
$V_{\iota+1}=\{0\}$, where $\iota$ is the step of $\G$.
For every $r>0$, a natural group automorphism $\delta_r:\cG\to\cG$
can be defined as the unique algebra homomorphism such that
$$
\delta_r(X):=r X\qquad\text{for every }X\in V_1.
$$
This one parameter group of mappings forms the family of the so-called
{\em dilations} of $\G$. Notice that simply connected nilpotent
Lie groups are diffeomorphic to their Lie algebra through the
exponential mapping $\exp:\cG\to\G$, hence dilations are automatically
defined as group isomorphisms of $\G$ and will be denoted by the
the same symbol $\delta_r$.

We will say that $\rho$ is a {\em homogeneous distance} on $\G$ if
it is a continuous distance of $\G$ satisfying the following conditions
\begin{equation}\label{defhomogmetr}
\rho(zx,zy)=\rho(x,y)\quad\text{and}\quad\rho(\delta_r(x),\delta_r(y))=r\rho(x,y)\qquad\text{for all }x,y,z\in\G,\:r>0.
\end{equation}
Important examples of homogeneous distances are the
well known Carnot-Carath\'eodory distance and the homogeneous
distance constructed in \cite{FSSC5}.

In the sequel, we will denote by $\haus^d$ and $\Shaus^d$, the $d$-dimensional Hausdorff and spherical Hausdorff measures
induced by a fixed homogeneous distance $\rho$, respectively.
Open balls of radius $r>0$ and centered at $x$ with respect to $\rho$
will be denoted by $B_{x,r}$ and the corresponding closed
balls will be denoted by $D_{x,r}$.
The number $Q$ denotes
the Hausdorff dimension of $\G$ with respect to $\rho$.

According to \eqref{stratificazione}, we say that an ordered set of vectors 
$$(X_1,X_2,\ldots,X_q)=(X^1_1,\dots X^1_{m_1},X^2_1,\dots,X^2_{m_2},\dots,X^\iota_1,\dots,X^\iota_{m_\iota})$$
 is an {\em adapted basis} of $\cG$ iff $m_k=$dim $v_k$ and
$$X^k_1,\dots,X^k_{m_k}$$
is a basis of the layer $V_k$ for every $k=1,\ldots,\iota$.
\begin{defi}
{\rm Let $(X_1,X_2,\ldots,X_q)$ be an adapted basis of $\cG$.
The {\em degree} $d(j)$ of $X_j$ is the unique integer $k$
such that $X_j\in V_k$.
Let 
$$
X_J:=X_{j_1}\wedge\dots\wedge X_{j_p}
$$
be a simple $p$-vector of $\Lambda_p\cG$, where
$J=(j_1,j_2,\ldots,j_p)$ and $1\leq j_1<j_2<\dots<j_p\leq q$.
The degree of $X_J$ is the integer $d(J)$ defined by the sum $d(j_1)+\dots+d(j_p)$.
}\end{defi}
Notice that the degree of a $p$-vector is independent from the
adapted basis we have chosen.
In the sequel, we will fix a {\em graded metric}
$g$ on $\G$, namely, a left invariant Riemannian metric on $\G$
such that the subspaces $V_k$'s are orthogonal.
It is easy to observe that all left invariant Riemannian metrics
such that $(X_1,\dots,X_q)$ is an orthonormal basis are graded metrics
and the family of $X_J$'s forms an orthonormal basis of
$\Lambda_p(\cG)$ with respect to the induced metric.
The norm induced by $g$ on $\Lambda_p(\cG)$ will be simply
denoted by $|\cdot|$.
\begin{defi}{\rm
When an adapted basis $(X_1,\ldots,X_q)$ is also orthonormal
with respect to the fixed graded metric $g$ is called {\em graded basis}.
}\end{defi}
\begin{defi}[Degree of $p$-vectors]{\rm
Let $\tau\in\Lambda_p(\cG)$ be a simple $p$-vector
and let $1\leq r\leq Q$ be a natural number.
Let $\tau=\sum_{J}\tau_J\;X_J$ be represented
with respect to the fixed adapted basis $(X_1,\ldots,X_q)$.
The projection of $\tau$ with degree $r$ is defined as
\begin{equation}\label{taur}
(\tau
)_r=\sum_{d(J)=r}\tau_J\;X_J.
\end{equation}
The {\em degree} of $\tau$ is defined as the integer
$$
d(\tau)=\max\left\{k\in\N\mid \mbox{such that $\tau_k\neq0$}\right\}.
$$
}\end{defi}
In the sequel, also an arbitrary auxiliary Riemannian metric
$\tilde g$ will be understood.
We define $\tau_\Sigma(x)$ as the unit
tangent $p$-vector to a $C^1$ submanifold 
$\Sigma$ at $x\in\Sigma$ with respect to the metric $\tilde g$, i.e.
$|\tau_\Sigma(x)|_{\tilde g}=1$. 
The {\em degree} of $x$ is defined as
\begin{eqnarray}\label{dsx}
	d_\Sigma(x)=d(\tau_\Sigma(x))
\end{eqnarray}
and the {\em degree of} $\Sigma$ is $d(\Sigma)=\max_{x\in\Sigma} d_\Sigma(x)$.
We will say that $x\in\Sigma$ has {\em maximum degree}
if $d_\Sigma(x)=d(\Sigma)$.
It is not difficult to check that these definitions
are independent from the fixed adapted basis $X_1,\dots,X_q$,
then they only depend on the tangent subbundle $T\Sigma$
and of the grading of $\cG$, namely they depend on
the ``geometric" position of the points with respect
to the grading \eqref{stratificazione}.
According to \eqref{taur}, we define ${\tau^d_{\Sigma}}(x)$ as the part of $\ts(x)$ with maximum degree $d=d(\Sigma)$, namely,
\begin{equation}\label{deftsm}
\tau_\Sigma^d(x)=\big(\tau_\Sigma(x)\big)_{d}.
\end{equation}
If $g$ is a fixed graded metric we will simply write
\begin{equation}\label{modulotsmax}
|\t_\Sigma^d(x)|=|\t_\Sigma^d(x)|_g.
\end{equation}
\begin{defi}\label{defpix}{\rm
Let $x\in\Sigma$ be a point of maximum degree.
Then we define
$$
\Pi_\Sigma(x)
=\{y\in\G:y=\exp(v)\text{ with }v\in\cG\text{ and }v\wedge\tau_\Sigma^d(x)=0\}\,.
$$
}\end{defi}
As a consequence of Corollary~\ref{vettangentemax},
we will see that $\Pi_\Sigma(x)$ is a subgroup of $\G$.
\subsection{Graded coordinates}
In the sequel the adapted basis $(X_1,\dots,X_q)$ will be fixed.
The exponential mapping $\exp:\cG\longrightarrow\G$ induces
a group law $C(X,Y)$ on $\cG$ for every $X,Y\in\cG$.
We have
\begin{eqnarray}\label{expXY}
\exp(X)\cdot\exp(Y)=\exp(C(X,Y)).
\end{eqnarray}
Recall that $C(X,Y)$ can be computed explicitly thanks to the Baker-Campbell-Hausdorff formula: for each multi-index of nonnegative integers $a=(a_1,\dots,a_l)$ we define
$$\begin{array}{l}
|a|:=a_1+\dots+a_l\\
\,a!:=a_1!\cdots a_l!
\end{array}$$
and we will say that $l$ is the length of $a$. If $b=(b_1,\dots,b_l)$ is another multi-index of length $l$ such that $a_l+b_l\geq 1$, and if $X,Y\in\cG$ we set
\begin{eqnarray*}
C_{ab}(X,Y) &:= & \left\{
\begin{array}{ll}
(\mbox{ad}\,X)^{a_1}(\mbox{ad}\,Y)^{b_1}\dots(\mbox{ad}\,X)^{a_l}(\mbox{ad}\,Y)^{b_l-1}\,Y & \mbox{if } b_l>0\\
(\mbox{ad}\,X)^{a_1}(\mbox{ad}\,Y)^{b_1}\dots(\mbox{ad}\,X)^{a_l-1}X
& \mbox{if } b_l=0.
\end{array}
\right.
\end{eqnarray*}
We used the notation $($ad$\,X)(Y):=[X,Y]$, agreeing that $($ad$\,X)^0$ is the identity.
According to \cite{Vara}), the Baker-Campbell-Hausdorff formula
is stated as follows
\begin{equation}\label{formulaBCH}
C(X,Y):=\sum_{l=1}^{\iota}\frac{(-1)^{l+1}}{l}\sum_{\substack{a=(a_1,\dots,a_l)\\
b=(b_1,\dots,b_l)\\ a_i+b_i\geq 1\:\forall i}} \frac{1}{a!b!|a+b|}C_{ab}(X,Y).
\end{equation}
For every adapted basis $(X_1,\dots,X_q)$, we can introduce a
system of {\em graded coordinates} on $\G$ given by
\begin{eqnarray}\label{Fco}
F:\R^q\longrightarrow\G,\quad F(x)=\exp\Big(\sum_{j=1}^qx_jX_j\Big),
\end{eqnarray}
where $\exp:\cG\longrightarrow\G$ is the exponential mapping.
Then the group law
\begin{eqnarray}\label{Fxy}
F(x)\cdot F(y)=F\big(P(x,y)\big)
\end{eqnarray}
is translated with respect to coordinates of $\R^q$ as 
\begin{equation}\label{defPeQ}
x\cdot y=P(x,y)=x+y+Q(x,y)\,,
\end{equation}
where the Baker-Campbell-Hausdorff formula \eqref{formulaBCH} implies that $P=(P_1,\dots,P_q)$ and $Q=(Q_1,\dots,Q_q)$ are polynomial vector fields.

It is also easy to check that dilations read in these coordinates as
$$\delta_r(x)=(rx_1,\dots,r^{d(j)}x_j,\dots,r^\iota x_q)\qquad\mbox{for every $r>0$}.$$
From definition of dilations and the Baker-Campbell-Hausdorff formula,
it follows that $Q_i(x,y)$ are homogeneous polynomials
with respect to dilations, i.e.
\begin{eqnarray}\label{homog}
P_i(\delta_r(x),\delta_r(y))=r^{d(i)}\,P_i(x,y)\quad\mbox{and}\quad
Q_i(\delta_r(x),\delta_r(y))=r^{d(i)}\,Q_i(x,y)\,.
\end{eqnarray}
As a result, we get
\begin{eqnarray}\label{qi}
\left\{\begin{array}{l}
Q_1=\dots=Q_{m_1}=0 \\
Q_i(x,y)=Q_i\left(\sum_{d(j)<i}x_j\;e_j,\sum_{d(j)<i}y_j\;e_j\right),
\end{array}\right.
\end{eqnarray}
where $(e_1,\ldots,e_q)$ denotes the canonical basis of $\R^q$
and $d(i)>1$.

Given a system of graded coordinates $F:\R^q\to\G$, we say that a function $p:\G\to\R$ is a {\em polynomial} on $\G$ if the composition $p\circ F^{-1}$ is a polynomial on $\R^q$; we say that $p$ is an {\em homogeneous polynomial of degree $l$} if it is a polynomial and $p(\delta_r(x))=r^l p(x)$ for any $x\in\G$ and $r>0$. It is not difficult to prove that $p$ is a homogeneous polynomial of degree $l$ if and only if $p\circ F^{-1}$ is a sum of monomials $$x_1^{l_1}x_2^{l_2}\cdots x_q^{l_q}\quad\text{with }\sum_{i=1}^q\,d(j) l_j=l.$$
Moreover, the notions of polynomial, homogeneous polynomial (and its degree) do not depend on the choice of graded coordinates $F$. Observe also that homogeneous polynomials of degree 0 are constants.

Any left invariant vector field $X_j$ of our fixed adapted basis
has a canonical representation as left invariant vector field $(F^{-1})_*(X_j)$ of $\R^q$, where $F$ is defined \eqref{Fco}.
We will use the same notation to indicate whis vector field in $\R^q$.
The left invariance of $X_j$ in $\R^q$ implies that
\begin{eqnarray*}
X_jf(x)= \partial_{y_j}(f\circ l_x)(0)
= df(x)\left( \frac{\partial P}{\partial y_j}(x,0)\right)\,,
\end{eqnarray*}
where $l_x(y)=x\cdot y\in\R^q$ and $f\in C^\infty(\R^q)$.
As a consequence, we have
\begin{equation}\label{formacampi}
X_j(x)=\sum_{i=1}^q X_{ij}(x)\;\partial_i=\displaystyle\sum_{i=1}^q \frac{\partial P_i}{\partial y_j}(x,0)\;\partial_i=\:\:\:\partial_j + \!\!\!\sum_{d(i)>d(j)} \frac{\partial Q_i}{\partial y_j}(x,0)\;\partial_i\,.
\end{equation}
By differentiating \eqref{homog} we get 
\begin{eqnarray}\label{homogvf}
X_{ij}(\delta_r(x))
=\frac{\partial P_i}{\partial y_j}(\delta_r(x),0)
=r^{d(i)-d(j)}\frac{\partial P_i}{\partial y_j}(x,0)
=r^{d(i)-d(j)}\,X_{ij}(x)\,,
\end{eqnarray}
i.e. $X_{ij}$ are homogeneous polynomials of degree $d(i)-d(j)$.

Next we present a key result in the proof of Lemma~\ref{curve}.
\begin{lem}\label{campisottoalg}
Let $J\subset\{1,2,\ldots,q\}$ be such that
$\mathcal F=\mbox{\rm span}\{X_j:j\in J\}$ is a subalgebra of $\cG$,
where $(X_1,\dots,X_q)$ be an adapted basis of $\cG$.
Then for every index $i\notin J$, the polynomial $Q_i(x,y)$ lies in the ideal generated by $\{x_l,y_l:l\notin J\}$, namely, we have
\begin{equation}\label{formacampisubalg}
Q_i(x,y)=\sum_{l\notin J,\: d(l)<d(i)} (x_lR_{il}(x,y) + y_lS_{il}(x,y))\,,
\end{equation}
where $R_{il},S_{il}$ are homogeneous polynomials of degree $d(i)-d(l)$.
\end{lem}
\begin{proof}
Let us fix $x,y\in\R^q$ and consider
$$
X:=\sum_{j=1}^q x_j X_j,\qquad Y:=\sum_{j=1}^q y_j X_j.
$$
By \eqref{expXY}, \eqref{Fxy}
and the Baker-Campbell-Hausdorff formula \eqref{formulaBCH},
we have
$$
C(X,Y)=\sum_{j=1}^q P_j(x,y)\;X_j\,,
$$
Therefore, defining $\pi_i:\cG\to\R$ as the
function which associates to every vector its $X_i$'s
coefficient, we clearly have $P_i(x,y)=\pi_i\big(C(X,Y)\big)$.
Thus, formulae \eqref{formulaBCH} and \eqref{defPeQ} yield
$$
Q_i(x,y)=\sum_{l=1}^{\iota}\frac{(-1)^{l+1}}{l}\sum_{\substack{a=(a_1,\dots,a_l)\\
b=(b_1,\dots,b_l)\\ a_i+b_i\geq 1\:\forall i}} \frac{1}{a!b!|a+b|} \pi_i(C_{ab}(X,Y)) - x_i - y_i.
$$
Observe that $C_{ab}(X,Y)$ is a commutator of $X$ and $Y$, whose length is equal to $|a+b|$; as the sum of commutator with length 1 gives $X+Y$ we get
$$
Q_i(x,y)=\sum_{l=1}^{\iota}\frac{(-1)^{l+1}}{l}\sum_{\substack{a=(a_1,\dots,a_l)\\
b=(b_1,\dots,b_l)\\ a_i+b_i\geq 1\:\forall i\\|a+b|\geq 2}} \frac{1}{a!b!|a+b|} \pi_i(C_{ab}(X,Y)).
$$
When the commutator $C_{ab}(X,Y)$ has length $h\geq 2$,
we can decompose it into the sum of commutators
of the vector fields $\{x_lX_l,y_lX_l : 1\leq l\leq q\}$.
Let us focus our attention on an individual addend of this sum
and consider its projection $\pi_i$. 
Clearly, this addend is a commutator of length $h$.
If this term is a commutator containing an element of the family
$\{x_lX_l,y_lX_l : l\notin J\}$, then its projection
$\pi_i$ will be a multiple of $x_l$ or $y_l$ for some $l\notin J$,
i.e. the projection $\pi_i$ of this term is a polynomial of the ideal
$$
\{x_l,y_l:l\notin J\}.
$$
On the other hand, if in the fixed commutator only elements of
$\{x_lX_l,y_lX_l : l\in J\}$ appear, then it belongs to $\mathcal F$.
In view of our hypothesis,
we have $\mathcal{F}\cap\mbox{\rm span}\{X_i\}=\{0\}$, hence
its projection through $\pi_i$ vanishes. This fact along with
(\ref{qi}) proves that $Q_i(x,y)$ has the form \eqref{formacampisubalg}.
\end{proof}
The next definition introduces the metric factor associated
with a simple $p$-vector. Notice that this definition generalized
the notion of metric factor first introduced in \cite{Mag2}.
\begin{defi}[Metric factor]\label{metrfct}{\rm
Let $\cG$ be a stratified Lie algebra equipped with a graded metric $g$
and a homogeneous distance $\rho$.
Let $\tau$ be a simple of $p$-vector of $\Lambda_p(\cG)$.
We define $\mathcal{L}(\tau)$ as the unique subspace associated with $\tau$.
The {\em metric factor} is defined by
\begin{eqnarray}\label{metricfactor}
\theta(\tau)=\mathcal{H}_{|\cdot|}^p
\left(F^{-1}\big(\exp\big(\mathcal{L}(\tau)\big)
\cap B_1\big)\right),
\end{eqnarray}
where $F:\R^q\longrightarrow\G$ is a system of graded coordinates with respect to an adapted orthonormal basis $(X_1,\dots,X_q)$.
The $p$-dimensional Hausdorff measure with respect to
the Euclidean norm of $\R^q$ has been denoted by 
$\mathcal{H}^p_{|\cdot|}$ and
$B_1$ is the open unit ball centered at $e$, with radius $r$
with respect to the fixed homogeneous distance $\rho$.
}\end{defi}

\section{Blow-up at points of maximum degree.}

\begin{lem}\label{buonivettoritangenti}
Let $\Sigma$ be a $p$-dimensional submanifold of class $C^1$
and let $x\in\Sigma$ be a point of maximum degree.
Then we can find
\begin{itemize}
\item a graded basis $X_1,\dots,X_q$ of $\cG$;
\item a neighbourhood $U$ of $x$;
\item a basis $v_1(y),\dots,v_p(y)$ of $T_y\Sigma$ for all $y\in U$
\end{itemize}
such that writing $v_j(y)=\sum_{i=1}^q C_{ij}(y) X_i(y)$, we have
\begin{equation}\label{matriceC}
C(y):=(C_{ij})_{\substack{i=1,\dots,q\\j=1,\dots,p}}=\left[\begin{array}{c|c|c|c}
Id_{\alpha_1} & 0 & \cdots & 0\\
O_1(y) & \ast & \cdots & \ast\\
\hline 0 & Id_{\alpha_2} & \cdots & 0 \\
0 & O_2(y) & \cdots & \ast\\
\hline \vdots & \vdots & \ddots & \vdots\\
\hline 0 & 0 & \cdots & Id_{\alpha_\iota}\\
0 & 0 & \cdots & O_\iota(y)
\end{array}\right]
\end{equation}
where $\alpha_k$ are integers satisfying
$0\leq\alpha_k\leq m_k$ and $\alpha_1+\dots+\alpha_\iota=p$.
The $(m_k-\alpha_k)\times\alpha_k$-matrix valued continuous functions $O_k$
vanish at $x$ and $\ast$ denotes a continuous bounded matrix
valued function.
\end{lem}
\begin{proof}
Observing that the degree of a point in $\Sigma$ is invariant
under left translations, it is not restrictive
assuming that $x$ coincides with the unit element $e$ of $\G$.

STEP 1. Here we wish to find the graded basis
$(X_1,\dots,X_q)$ of $\cG$ and the basis $v_1,\dots,v_p$ of $T_e\Sigma$
required in the statement of the lemma and that
satisfy \eqref{matriceC} when $y=e$.
Let us fix a basis $(t_1,\dots,t_p)$ of $T_e\Sigma$ and
use the same notation to denote the corresponding basis
of left invariant vector fields of $\cG$.
We denote by $\pi_k$ the canonical projection of
$\cG$ onto $V_k$. Let $0\leq\alpha_\iota\leq m_\iota$ be the dimension of the subspace spanned by 
$$\pi_\iota(t_1),\dots,\pi_\iota(t_q).$$
Taking linear combinations of $t_j$ we can suppose that
the first $\alpha_\iota$ vectors $\{\pi_\iota(t_j)\}_{1\leq j \leq \alpha_\iota}$ form an orthonormal set of $V_\iota$,
with respect to the fixed graded metric $g$.
Then we set
$$
X_j^\iota:=\pi_\iota(t_j)\in V_\iota\quad\mbox{and}\quad
v_j^\iota:=t_j\in T_e\Sigma\,,
$$
whenever $1\leq j \leq \alpha_\iota$.
Adding proper linear combinations of these $t_j$ to the remaining 
vectors of the basis, we can assume that
$\{t^{\iota-1}_j:=t_{j+\alpha_\iota}\}_{1\leq j \leq p-\alpha_\iota}$ are linearly independent and that
$$
\pi_\iota(t^{\iota-1}_j)=0\qquad \mbox{whenever}\quad
j=1,\dots,p-\alpha_\iota.
$$
Now consider the $p-\alpha_\iota$ vectors
$$
\pi_{\iota-1}(t^{\iota-1}_1),\dots,\pi_{\iota-1}(t^{\iota-1}_{p-\alpha_\iota})
$$
and let $0\leq\alpha_{\iota-1}\leq m_{\iota-1}$ be the rank of the subspace of $V_{\iota-1}$ generated by these vectors.
Taking linear combinations of $t^{\iota-1}_j$, we can suppose that
$\pi_{\iota-1}(t^{\iota-1}_j)$ with $j=1,\ldots,\alpha_{\iota-1}$
form an orthonormal set of $V_{\iota-1}$ and that defining
$
\{t^{\iota-2}_j:=t^{\iota-1}_{j+\alpha_{\iota-1}}\}
_{1\leq j \leq p-\alpha_\iota-\alpha_{\iota-1}}
$
we have
$$
\pi_{\iota-1}(t^{\iota-2}_j)=0\quad\mbox{whenever}\quad j=1,\dots,p-\alpha_\iota-\alpha_{\iota-1}.
$$
Then we set
$$
X_j^{\iota-1}:=\pi_{\iota-1}(t^{\iota-1}_j)\in V_{\iota-1}
\quad\mbox{and}\quad v_j^{\iota-1}:=t^{\iota-1}_j\in T_e\Sigma\,.
$$
for every $j=1,\ldots,\alpha_{\iota-1}$.
Repeating this argument in analogous way, we obtain integers $\alpha_k$
with $0\leq\alpha_k\leq m_k$ for every $k=1,\ldots,\iota$ and
vectors
$$
X_j^k\in V_k\,,\qquad v_j^k\in T_e\Sigma,
\quad\mbox{where}\quad k=1,\dots,\iota\quad\mbox{and}\quad
\:j=1,\dots,\alpha_k.
$$
Notice that $\alpha_1+\dots+\alpha_\iota=p$ and that 
\begin{eqnarray}\label{tangentsigma}
(v_1^1,\dots,v_{\alpha_1}^1,\dots,v_1^\iota,\dots,v_{\alpha_\iota}^\iota)
\end{eqnarray}
is a basis of $T_e\Sigma$. 
We complete the $X_j^k$'s to a graded basis 
$$
(X^1_1,\dots X^1_{m_1},X^2_1,\dots,X^2_{m_2},
\dots,X^\iota_1,\dots,X^\iota_{m_\iota})
$$
of $\cG$, that will be also denoted by $(X_1,\dots,X_q)$.
It is convenient to relabel the basis \eqref{tangentsigma}
as  $(v_1,\dots,v_p)$, hence we write
$v_j=\sum_{i=1}^q C_{ij} X_i$ obtaining
$$C:=(C_{ij})=\left[\begin{array}{c|c|c|c}
Id_{\alpha_1} & \ast & \cdots & \ast\\
0 & \ast & \cdots & \ast\\
\hline 0 & Id_{\alpha_2} & \cdots & \ast \\
0 & 0 & \cdots & \ast\\
\hline \vdots & \vdots & \ddots & \vdots\\
\hline 0 & 0 & \cdots & Id_{\alpha_\iota}\\
0 & 0 & \cdots & 0
\end{array}\right].$$
Performing suitable linear combinations of $v_j$'s, we can
assume that 
\begin{equation}\label{Cin0}
C=\left[\begin{array}{c|c|c|c}
Id_{\alpha_1} & 0 & \cdots & 0\\
0 & \ast & \cdots & \ast\\
\hline 0 & Id_{\alpha_2} & \cdots & 0 \\
0 & 0 & \cdots & \ast\\
\hline \vdots & \vdots & \ddots & \vdots\\
\hline 0 & 0 & \cdots & Id_{\alpha_\iota}\\
0 & 0 & \cdots & 0
\end{array}\right].
\end{equation}

STEP 2. The basis $(v_1,\ldots,v_p)$ of $T_e\Sigma$
can be extended to a frame of continuous vector fields
$(v_1(y),\dots,v_p(y))$ on $\Sigma$ defined in neighborhood $U$ of $e$.
Thanks to the previous step, defining $v_j(y)=\sum_{i=1}^q C_{ij}(y) X_i(y)$
we have
$$
C(y):=(C_{ij}(y))=\left[\begin{array}{c|c|c|c}
Id_{\alpha_1} + o(1) & o(1) & \cdots & o(1)\\
o(1) & \ast & \cdots & \ast\\
\hline o(1) & Id_{\alpha_2} +o(1) & \cdots & o(1) \\
o(1) & o(1) & \cdots & \ast\\
\hline \vdots & \vdots & \ddots & \vdots\\
\hline o(1) & o(1) & \cdots & Id_{\alpha_\iota} + o(1)\\
o(1) & o(1) & \cdots & o(1)
\end{array}\right]
$$
where $o(1)$ denotes a matrix-valued continuous function vanishing at $e$.
Observing that $Id_{\alpha_k}+o(1)$ are still invertible for every
$y$ in a smaller neighbourhood $U'\subset U$ of $e$,
we can replace the $v_j$'s with linear combinations to get
$$C(y)=\left[\begin{array}{c|c|c|c}
Id_{\alpha_1}+ o(1) & 0 & \cdots & 0\\
o(1) & \ast & \cdots & \ast\\
\hline
 0 & Id_{\alpha_2}+ o(1) & \cdots & 0 \\
o(1) & o(1) & \cdots & \ast\\
\hline
\vdots & \vdots & \ddots & \vdots\\
\hline
0 & 0 & \cdots & Id_{\alpha_\iota}+ o(1) \\
o(1) & o(1) & \cdots & o(1)
\end{array}\right].$$
The same argument leads us to define a new frame with matrix
\begin{equation}\label{Cfinale}
C(y)=\left[\begin{array}{c|c|c|c}
Id_{\alpha_1} & 0 & \cdots & 0\\
O_1(y) & \ast & \cdots & \ast\\
\hline
0 &Id_{\alpha_2}& \cdots & 0 \\
o(1) & O_2(y) & \cdots & \ast\\
\hline
\vdots & \vdots & \ddots & \vdots\\
\hline
0 & 0   & \cdots & Id_{\alpha_\iota} \\
o(1) & o(1)   & \cdots & O_\iota(y)
\end{array}\right]\,,
\end{equation}
where $O_j$ are defined in the statement of the present lemma.
To finish the proof, it remains to show
that all $o(1)$'s of (\ref{Cfinale})
are actually null matrix functions.
Here we utilize the fact that the submanifold has
maximum degree at $e$.
Notice that the simple $p$-vector
$$
v_1(y)\wedge\dots\wedge v_p(y)=\sum_J a_J(y) X_J(y)\,,
$$
is proportional to the tangent vector $\tau_\Sigma(y)$.
In addition, if $J=(j_1,\dots,j_p)$, then
$a_J(y)$ is the determinant of the $p\times p$ submatrix
obtained taking the $j_1$-th, $j_2$-th, $\dots,j_{p-1}$-th
and $j_p$-th row of $C(y)$.
From \eqref{Cin0} we immediately conclude that $d_\Sigma(e)=\alpha_1+2\alpha_2+\dots+\iota\alpha_\iota$.
Finally, where one entry of some $o(1)$
does not vanish, it is possible to find some $J_0$
such that $d(J_0)>\alpha_1+2\alpha_2+\dots+\iota\alpha_\iota$ and $a_J(y)\neq 0$. This would imply $d_\Sigma(y)>d_\Sigma(e)$, contradicting
the assumption that $d_\Sigma(e)=\max_{y\in U'}d_\Sigma(y)$.
\end{proof}
\begin{rem}
{\upshape It is easy to interpret the statement and the proof of
Lemma~\ref{buonivettoritangenti} in the case some $\alpha_k$
vanishes. Clearly, the $\alpha_k$ columns in \eqref{matriceC} intersecting
$I_{\alpha_k}$ and then the corresponding vectors $v^k_j$
disappear.}
\end{rem}
\begin{rem}\label{piuregolarita}
{\upshape Clearly, when $\Sigma$ is of class $C^r$
the $v_j$'s of the previous lemma are of class $C^{r-1}$.
In fact, the linear transformations performed in
the proof of Lemma~ \ref{buonivettoritangenti}
are of class $C^{r-1}$.}
\end{rem}
The previous lemma allows us to state the following definitions.
\begin{defi}{\rm
Let $\Sigma$ be a $C^1$ smooth submanifold and let
$x\in\Sigma$ be a point of maximum degree.
Then we can define the {\em degree} $\sigma:\{1,\ldots,p\}\longrightarrow\N$
induced by $\Sigma$ at $x$ as follows
$$
\sigma(j)=i\qquad\mbox{if}\qquad \sum_{s=1}^{i-1}\alpha_s< j
\leq\sum_{s=1}^i\alpha_s,
$$
where $\alpha_i$ are defined in Lemma~\ref{buonivettoritangenti}.
}\end{defi}
\begin{defi}{\rm
Let $\Sigma$ be a $C^1$ smooth submanifold and let
$x\in\Sigma$ be a point of maximum degree.
Then we will denote by 
$$
(X_1^1,\dots,X^1_{m_1},\dots,X^\iota_1,\dots,X_{m_\iota}^\iota)
\quad\mbox{and}
\quad (v_1^1,\dots,v^1_{\alpha_1},\dots,v^\iota_1,\dots,v_{\alpha_\iota}^\iota)
$$
the frames on $\G$ and on a neighborhood $U$ of $z$ in $\Sigma$, respectively, which satisfy the conditions
of Lemma~\ref{buonivettoritangenti}.
We will also indicate these frames by 
$$
(X_1,\dots,X_q)\qquad\mbox{and}\qquad (v_1,\dots,v_p).
$$
}\end{defi}
\begin{cor}\label{vettangentemax}
Let $\Sigma$ be a $C^1$ smooth submanifold with
$x\in\Sigma$ satisfying $d_\Sigma(x)=d(\Sigma)$.
Then  ${\tau_{\Sigma}^d}(x)$ is a simple
$p$-vector which is proportional to
$$
X_1^1\wedge\dots\wedge X_{\alpha_1}^1\wedge\dots\wedge X^\iota_1\wedge\dots\wedge X_{\alpha_\iota}^\iota\,,
$$
then we also have
$$
\Pi_\Sigma(x)=\exp\bigl(\mbox{\rm span}
\{X_1^1,\dots,X^1_{\alpha_1},\dots,X^\iota_1,\dots,X_{\alpha_\iota}^\iota\}
\bigr)
$$
\end{cor}
\begin{proof}
By expression \eqref{matriceC}, $\tau_\Sigma$ is clearly
proportional to 
\begin{equation}
X_1^1\wedge\dots\wedge X_{\alpha_1}^1\wedge\dots\wedge X^\iota_1\wedge\dots\wedge X_{\alpha_\iota}^\iota + R\,,
\end{equation}
where $R$ is a linear combination of simple $p$-vectors
with degree less than $d(X_1^1\wedge\dots\wedge X_{\alpha_\iota}^\iota)$. Then $d=d(X_1^1\wedge\dots\wedge X_{\alpha_\iota}^\iota)$ and
${\tau_{\Sigma}^d}(x)$ is proportional to
$X_1^1\wedge\dots\wedge X_{\alpha_\iota}^\iota$.
\end{proof}
\begin{defi}{\rm
We will denote by
\begin{eqnarray}\label{basispi}
(X_1^1,\dots,X^1_{\alpha_1},\dots,X^\iota_1,\dots,X_{\alpha_\iota}^\iota)
\end{eqnarray}
the frame of Corollary~\ref{vettangentemax}, arising from Lemma~\ref{buonivettoritangenti}, and by
\begin{eqnarray}\label{canproj}
\pi_\Sigma(x):\G\longrightarrow\Pi_\Sigma(x)
\end{eqnarray}
the corresponding canonical projection.
}\end{defi}
\begin{cor}\label{graficosulpiano}
Let $e\in\Sigma$ be such that $d_\Sigma(e)=d(\Sigma)$.
Let us embed $\Sigma$ into $\R^q$ by the system of graded
coordinates $F$ induced by $\{X_j^k\}_{k=1,\dots,\iota,\:j=1,\dots,m_k}$.
Then there exists a function 
\begin{eqnarray*}
\f &:& A\subset\R^p\longrightarrow\R^{q-p}\\
&&x=(x_1^1,\dots,x_{\alpha_1}^1,\dots,x_{\alpha_\iota}^\iota)
\longmapsto (\f_{\alpha_1+1}^1,\dots,\f_{m_1}^1,\dots,\
\f_{\alpha_\iota+1}^\iota,\dots,\f_{m_\iota}^\iota)(x),
\end{eqnarray*}
defined on an open neighbourhood $A\subset\R^p$ of zero,
such that $\f(0)=0$ and $\Sigma\supset\Phi(A)$,
where $\Phi$ is the mapping defined by
\begin{eqnarray}\label{repinv}
&&\F :A\to\R^q \nonumber\\
&& x\longrightarrow \big(x_1^1,\dots,x_{\alpha_1}^1,\f_{\alpha_1+1}^1(x),\dots,\f_{m_1}^1(x),\dots,x_1^\iota,\dots,x_{\alpha_\iota}^\iota,\f_{\alpha_\iota+1}^\iota(x),\dots,\f_{m_\iota}^\iota(x)\big).
\end{eqnarray}
and satisfying $\nabla\F(0)=C(0)$, with $C$ given by
Lemma \ref{buonivettoritangenti}.
\end{cor}
\begin{proof}
Representing $\pi_\Sigma(x)$ with respect to our graded coordinates,
we obtain
\begin{eqnarray*}
\tilde\pi_\Sigma(x) &:& \R^q\to\R^p\\
&& x\longmapsto (x_1^1,\dots,x_{\alpha_1}^1,\dots,x_1^\iota,\dots,x_{\alpha_\iota}^\iota)\,.
\end{eqnarray*}
Taking its restriction 
\begin{eqnarray*}
\pi &:& \Sigma\to\R^p\\
&& x\longmapsto (x_1^1,\dots,x_{\alpha_1}^1,\dots,x_1^\iota,\dots,x_{\alpha_\iota}^\iota)\,,
\end{eqnarray*}
we wish to prove that $\pi$
is invertible near $0$, i.e.
that $d\pi(0):T_0\Sigma\to\R^p$ is onto.
According to \eqref{matriceC} and the fact that
$\pi$ is the restriction of a linear mapping, it follows that
$d\pi(v_j^k(0))=\partial_{x_j^k}$
for every $k=1,\dots,\iota$ and $j=1,\dots,\alpha_k$.
This implies the existence of $\Phi=\pi_{|U}^{-1}$ having the
representation \eqref{repinv}, hence one can easily check that
$d\pi(\partial_{x_j^k}\!\Phi\,(0))=\partial_{x_j^k}$ also holds
for every $k=1,\dots,\iota$ and $j=1,\dots,\alpha_k$.
As a consequence, invertibility of $d\pi(0):T_0\Sigma\rightarrow\R^p$ 
gives $v_j^k(0)=\partial_{x_j^k}\!\Phi\,(0)$.
It follows that each column of $\nabla\F(0)$ equals the
corresponding one of $C(0)$.
\end{proof}
From now on, we will assume that $\Sigma$ is a $C^{1,1}$ submanifold
of $\mathbb G$.
\begin{lem}\label{subalgebra}
Let $x\in\Sigma$ be such that $d_\Sigma(x)=d(\Sigma)$.
Then $\Pi_\Sigma(x)$ is a subgroup.
\end{lem}
\begin{proof}
Posing $d=d(\Sigma)$, due to Corollary~\ref{vettangentemax}, ${\tau_{\Sigma}^d}(x)$
is proportional to the simple $p$-vector
$$
X_1^1\wedge\dots\wedge X^1_{\alpha_1}\wedge\cdots
\wedge X^\iota_1\wedge\cdots\wedge X_{\alpha_\iota}^\iota\,.
$$
We define $\mathcal F$ as the space of linear combinations of
vectors $\{X^k_j\}^{k=1,\ldots,\iota}_{j=1,\ldots,\alpha_k}$.
It suffices to prove that each bracket
$[X_j^k,X_i^l]$ lies in $\mathcal F$ for every $1\leq k,l \leq \iota$,
$1\leq j\leq \alpha_k$ and $1\leq i\leq \alpha_l$.
Taking into account Remark \ref{piuregolarita},
we can find Lipschitz functions $\phi_r,\psi_s$, which vanish at $x$ whenever $d(r)=k$ or $d(s)=l$, such that
$$
v_j^k=X_j^k+\sum_{d(r)\leq k} \phi_r\;X_r\quad\mbox{and}\quad
v_i^l=X_i^l+\sum_{d(s)\leq l} \psi_s\;X_s.
$$
For a.e. $y$ belonging to a neighbourhood $U$ of $x$,
we have 
\begin{eqnarray}
[v_j^k,v_i^l]&=&\Big[X_j^k+\sum_{d(r)\leq k} \phi_r\; X_r,\:\:
X_i^l+\sum_{d(s)\leq l} \psi_s\; X_s\Big]\nonumber\\
&=& [X_j^k,X_i^l] +\sum_{d(r)\leq k}\phi_r\;[X_r,X^l_i]+
\sum_{d(s)\leq l}\psi_s\;[X^k_j,X_s] + \sum_{d(r)\leq k, d(s)\leq l}
\phi_r\;\psi_s\;[X_r,X_s] \nonumber \\
&+&\sum_{d(s)\leq l} (X_j^k\psi_s) X_s
- \sum_{d(r)\leq k} (X_i^l\phi_r) X_r \\
&+& \sum_{d(r)\leq k, d(s)\leq l} \Big(\phi_r\,(X_r\psi_s)\; 
X_s\;-\;\psi_s\,(X_s\psi_r)\;X_r\Big) \nonumber
\end{eqnarray}
By Frobenius theorem we know that this vector is tangent to $\Sigma$, i.e. it is a linear combination of $v_1^1,\dots,v^\iota_{\alpha_\iota}$ and
lies in $V_1\oplus\dots\oplus V_{k+l}$, hence
Lemma~\ref{buonivettoritangenti} implies that it must be of the form
$$
[v_j^k,v_i^l] = \sum_{\sigma(r)\leq k+l} a_r v_r.
$$
Projecting both sides of the previous identity onto $V_{k+l}$, we get
\begin{multline*}
[X_j^k,X_i^l] +\sum_{d(r)= k}\phi_r\;[X_r,X^l_i]+
\sum_{d(s)= l}\psi_s\;[X^k_j,X_s] +\\+ \sum_{d(r)= k, d(s)= l}
\phi_r\;\psi_s\;[X_r,X_s]=\sum_{\sigma(r)= k+l} a_r\; \pi_{k+l}(v_r).
\end{multline*}
From \eqref{matriceC} the projections
$\pi_{k+l}\big(v_r(y)\big)$ converge to a linear combination
of vectors $X_i^{k+l}$ as $y$ goes to $x$, where $1\leq i\leq\alpha_{k+l}$.
We can find a sequence of points
$(y_\nu)$ contained in $U$, where $[v_j^k,v_i^l]$ is defined and $y_\nu\rightarrow x$ as $\nu\rightarrow\infty$. 
Then the coefficients $a_r$ are defined on $y_\nu$ and up to
extracting subsequences it is not restrictive assuming that
$a_r(y_\nu)$, which is bounded since $\Sigma$ is $C^{1,1}$, converges for every $r$ such that $\sigma(r)\leq k+l$.
Thus, restricting the previous equality on the set $\{y_\nu\}$
and taking the limit as $\nu\rightarrow\infty$, it follows that
$[X_j^k,X_i^l]$ is a linear combination of $\{X_i^{k+l}\}_{1\leq i\leq\alpha_{k+l}}$. This ends the proof. 
\end{proof}
Let us consider the parameters $\lambda=(\lambda_1^1,\dots,\lambda_{\alpha_1}^1,
\dots,\lambda_1^\iota,\dots,\lambda_{\alpha_\iota}^\iota)\in\R^p$
and a point $e\in\Sigma$ with $d_\Sigma(e)=d(\Sigma)$.
We aim to study properties of solution
$\gamma(t,\lambda)$ of the Cauchy problem
\begin{equation}\label{gammalambda}
\left\{\begin{array}{l}
\partial_t\gamma(t,\lambda)= \displaystyle\sum_{\substack{k=1,\dots,\iota\\j=1,\dots,\alpha_k}} \lambda_j^k\: v_j^k\big(\gamma(t,\lambda)\big)\, t^{k-1}\\
\gamma(0,\lambda)=0
\end{array}\right.\,,
\end{equation}
where the vector fields $v^k_j$ are defined in Lemma~\ref{buonivettoritangenti} with $x=e$.

For every compact set $L\subset \R^p$, there exists
a positive number $t_0=t_0(L)$ such that $\g(\cdot,\lambda)$
is defined on $[0,t_0]$ for every $\lambda\in L$. 

The next lemma gives crucial estimates on
the coordinates of $\gamma(\cdot,\lambda)$.
Notice that graded coordinates arising from the corresponding
graded basis $(X_1,\ldots,X_q)$ will be understood.

\begin{lem}\label{curve}
Let $\gamma(\cdot,\lambda)$ be the solution of (\ref{gammalambda}).
Then for every $k=1,\dots,\iota$ and every $j=1,\dots,m_k$
there exist homogeneous polynomials $g_j^k$ of degree $k$,
that vanish when $k=1$, have the form
$g_j^k(\lambda_1^1,\dots,\lambda^1_{\alpha_1},\ldots,
\lambda^{k-1}_1,\ldots,\lambda_{\alpha_{k-1}}^{k-1})$
when $k>1$, satisfy $g_j^k(0)=0$ and the estimates
\begin{equation}\label{stimesopra}
\gamma_j^k(t,\lambda)=\left\{\begin{array}{ll}
\Big(\l_j^k/k + g_j^k(\l_1^1,\dots,\l_{\alpha_{k-1}}^{k-1})\Big)\; t^k + O(t^{k+1}) & \text{if }1\leq j\leq \alpha_k\vspace{.1cm}\\
O(t^{k+1}) & \text{if }\alpha_k+1\leq j\leq m_k
\end{array}\right. 
\end{equation}
hold for every $\lambda\in L$ and every $t\in[0,t_0]$.
\end{lem}
\begin{proof}
From \eqref{formacampi} and \eqref{homogvf},
we have $X_s=\sum_{i=1}^q X_{is}\,e_i$ where
\begin{eqnarray}\label{bpro}
X_{is}(x)=\left\{\begin{array}{lcc}
\delta_{is} & \mbox{if} & d(i)\leq d(s) \\
u_{is}(x^1_1,\ldots,x^1_{m_1},\ldots,x^{d(i)-1}_1
\!\!\!\!\!\!\!\!,\ldots,x^{d(i)-1}_{m_{d(i)}-1})   & \mbox{if} & d(i)>d(s) 
\end{array}\right.
\end{eqnarray}
and $u_{is}$ is a homogeneous polynomial satisfying
$u_{is}(\delta_r(x))=r^{d(i)-d(s)}u_{is}(x)$.
Setting
$$
\begin{array}{l}
\tilde\lambda=\tilde\l(t)=(\lambda_1^1,\dots,\l_{\alpha_1}^1,\l_1^2 t,\dots,\l_{\alpha_2}^2 t, \dots,\l_1^\iota t^{\iota-1},\dots,\l_{\alpha_\iota}^\iota t^{\iota-1}) \in\R^p
\end{array}
$$
and taking into account the expression of $v_j$ given in
Lemma~\ref{buonivettoritangenti}, we can write the Cauchy problem (\ref{gammalambda}) as 
\begin{equation}\label{gammapunto}
\partial_t\gamma(t,\lambda)=\sum_{r=1}^p v_r\big(\gamma(t,\lambda)\big)\,\tilde\lambda_r(t)=
\sum_{r=1}^p \sum_{s=1}^q C_{sr}\big(\gamma(t,\lambda)\big)\,
X_s\big(\gamma(t,\lambda)\big) \tilde\l_r(t)\,,
\end{equation}
where $C(\cdot)$ is given by Lemma \ref{buonivettoritangenti}.
Now we fix $\l\in L$ and write for simplicity $\gamma$
in place of $\g(\cdot,\lambda)$. The coordinates of $\gamma$
will be also denoted as follows
$$
(\g_1^1,\dots,\g_{m_1}^1,\dots,\gamma^\iota_1,\ldots,\gamma_{m_\iota}^\iota).
$$

STEP 1. We start proving \eqref{stimesopra} for the coordinates of $\gamma$
belonging to the first layer, i.e. 
\begin{equation}\label{stimestrato1}
\left\{\begin{array}{ll}
\g_j^1(t) = \l_j^1 t & \text{if }1\leq j\leq \alpha_1\vspace{.1cm}\\
\g_j^1(t) = O(t^2) & \text{if }\alpha_1+1\leq j\leq m_1
\end{array}\right.\,.
\end{equation}
In view of \eqref{gammapunto}, we get 
$$
\dot\g_j^1 = \sum_{r=1}^p \sum_{s=1}^q C_{sr}(\g)\,X_{js}(\g)\,
\tilde\l_r.
$$
For $1\leq j\leq \alpha_1$ we have
$1=d(j)\leq d(s)$, then \eqref{bpro} imply that
$X_{js}=\delta_{js}$, whence  
$$
\dot\g_j^1 = \sum_{r=1}^p C_{jr}(\g) \tilde\l_r = \tilde\l_j=\l_j^1,
$$
where the second equality follows from \eqref{matriceC}, which
implies $C_{jr}(x)=\delta_{jr}$.
This shows the first equality of \eqref{stimestrato1}.
Now we consider the case $\alpha_1+1\leq j\leq m_1$.
Due to \eqref{bpro} and $1=d(j)\leq d(s)$, we have
\begin{eqnarray}\label{twoterms}
\dot\g_j^1 = \sum_{r=1}^p C_{jr}(\gamma) \tilde\l_r
= \sum_{\sigma(r)=1} C_{jr}(\g) \tilde\l_r
+\sum_{\sigma(r)\geq 2} C_{jr}(\g) \tilde\lambda_r\,.
\end{eqnarray}
From \eqref{matriceC}, we have $C_{jr}(y)=o(1)$ whenever $\sigma(r)=1$, hence $C_{jr}(\gamma(t))=o(t)$.
From the same formula, we deduce that $C_{jr}(x)$ is bounded
whenever $\sigma(r)\geq 2$, and for the same indices $r$
we also have $\tilde\lambda_r=O(t)$, hence the second adddend
of (\ref{twoterms}) is equal to $O(t)$.
We have shown that $\dot\gamma_j^1=O(t)$ for every
$\alpha_1+1\leq j\leq m$, therefore the second equality
of \eqref{stimestrato1} is proved.

STEP 2. We will prove \eqref{stimesopra} by induction on
$k=1,\dots,\iota$. The previous step yields these estimates
for $k=1$. Let us fix $k\geq 2$ and suppose that \eqref{stimesopra}
holds for all integers less than or equal to $k-1$.
Next, we wish to prove \eqref{stimesopra} for components of
$\gamma$ with degree $k$ and for any fixed $1\leq j\leq m_k$.
We denote by $i$ the unique integer between 1 and $q$ such that $X_i=X_j^k$
and accordingly we have $\gamma_i=\gamma_j^k$, where $d(i)=k$.
Taking into account (\ref{bpro}) and that
$C_{sr}$ vanishes when $d(s)>\sigma(r)$, it follows that
\begin{eqnarray}\label{gsum}
\dot\gamma_i &=& \sum_{r=1}^p \sum_{s=1}^q X_{is}(\g) C_{sr}(\g) \tilde\l_r
= \sum_{\substack{1\leq r\leq p\\ d(s)\leq d(i)\\ d(s)\leq \sigma(r)}} X_{is}(\g) C_{sr}(\g) \tilde\l_r\,.
\end{eqnarray}
We split this sum into three addends
\begin{equation}\label{3split}
\dot\gamma_j^k=\dot\gamma_i
= \sum_{\substack{1\leq r\leq p\\ d(i)\leq\sigma(r)}} C_{ir}(\g) \tilde\l_r +
\sum_{\substack{1\leq r\leq p\\ d(s)<d(i)\\ d(s)=\sigma(r)}} X_{is}(\g) C_{sr}(\g) \tilde\l_r +
\sum_{\substack{1\leq r\leq p\\ d(s)< d(i)\\ d(s)<\sigma(r)}}
X_{is}(\g) C_{sr}(\g) \tilde\l_r\,.
\end{equation}
We first consider the case $1\leq j\leq\alpha_k$.
Then \eqref{matriceC} implies that $C_{ir}(x)=\delta_{ir}$,
therefore the first term of \eqref{3split}
coincides with $\tilde\lambda_i(t)=\l_j^k t^{k-1}$.
Now we deal with the remaining terms.
Our inductive hypothesis yields
\begin{eqnarray}\label{sopraind}
&&\gamma_s^l(t,\lambda)=\left\{\begin{array}{ll}
\bigl(\l_s^l/l + g_s^l(\l_1^1,\dots,\lambda_{\alpha_{l-1}}^{l-1})
+ O(t)\bigr)\;t^l & \text{if }1\leq s\leq \alpha_l\vspace{.1cm}\\
O(t)\;t^l & \text{if }\alpha_l+1\leq s\leq m_l
\end{array}\right.\,, 
\end{eqnarray}
whenever $l\leq k-1$, where $g^l_s$ is a homogeneous polynomial
of degree $l$. Due to (\ref{bpro}), $X_{is}$ are homogeneous polynomials
of degree $d(i)-d(s)=k-d(s)>0$, then applying \eqref{sopraind}, we achieve
\begin{eqnarray}\label{hompolprof}
X_{is}(\gamma_1^1,\dots,\gamma_{m_{k-1}}^{k-1})
= \bigl(N_{is}(\l_1^1,\dots,\l_{\alpha_{k-1}}^{k-1})
+O(t)\bigr)\;t^{k-d(s)}
\end{eqnarray}
whenever $d(s)\leq d(i)=k$ and $\upsilon_{is}=\delta_{is}$ if $d(s)=k$.
Notice that $N_{is}$ are homogeneous polynomial of degree $k-d(s)$
since it is a composition of the homogeneous polynomial
$X_{is}$ and of the homogeneous polynomials $\l_s^l/l + g_s^l(\l_1^1,\dots,\lambda_{\alpha_{l-1}}^{l-1})$ with degree $l$.

Let us focus our attention on the second addend of \eqref{3split}.
By definition of $\tilde\lambda$, we have
$\tilde\l_r= \l^{\sigma(r)}_{l(r)} t^{\sigma(r)-1}$,
for some $1\leq l(r)\leq \alpha_{\sigma(r)}$, hence
this second term equals 
\begin{eqnarray*}
&&\sum_{\substack{1\leq r\leq p\\ d(s)< d(i)\\ d(s)= \sigma(r)}}
\bigl[C_{sr}(0)+O(t)\bigr] \bigl[N_{is}(\lambda_1^1,\dots
\l_{\alpha_{k-1}}^{k-1}) t^{k-d(s)}
+O(t^{k-d(s)+1})\bigr] \lambda_{l(r)}^{\sigma(r)} t^{\sigma(r)-1}\\
&&=\sum_{\substack{1\leq r\leq p\\ d(s)< d(i)\\ d(s)= \sigma(r)}}
C_{sr}(0)\, N_{is}(\lambda_1^1,\dots\lambda_{\alpha_{k-1}}^{k-1})
\, \lambda_{l(r)}^{d(s)}\, t^{k-1}+O(t^k)
\\
&&= \tilde N_i(\l_1^1,\dots,\l_{\alpha_{k-1}}^{k-1})\, t^{k-1} + O(t^k),
\end{eqnarray*}
where $\tilde N_i$ is a homogeneous polynomial of degree $k=d(i)$.
From \eqref{hompolprof} and taking into account the definition of $\tilde\lambda_r$, the last term of \eqref{3split} can be
written as follows
\begin{eqnarray*}
&& \sum_{\substack{1\leq r\leq p\\ d(s)< d(i)\\ d(s)<\sigma(r)}}
C_{sr}(\gamma(t)) \bigl[N_{is}(\l_1^1,\dots,\l_{\alpha_{k-1}}^{k-1})
t^{k-d(s)}+O(t^{k-d(s)+1})\bigr]  O(t^{\sigma(r)-1})\\
&&= \sum_{\substack{1\leq r\leq p\\ d(s)< d(i)\\ d(s)<\sigma(r)}}
O(t^{k-d(s)+\sigma(r)-1})=O(t^k)\,.
\end{eqnarray*}
Summing up the results obtained for the three addends of \eqref{3split},
we have shown that
$$
\dot\g_j^k(t) = (\l_j^k + \tilde N_i(\lambda_1^1,\dots,
\lambda_{\alpha_{k-1}}^{k-1}))t^{k-1} + O(t^k)
$$
whence the first part of \eqref{stimesopra} follows. 

Next, we consider the case $\alpha_k+1\leq j\leq m_k$.
In this case we decompose (\ref{gsum}) into the following
two addends
\begin{eqnarray}\label{2split}
\dot\gamma_i &=&
\sum_{\substack{1\leq r\leq p\\ k\leq \sigma(r)}}
C_{ir}(\g)\; \tilde\l_r\,+
\sum_{\substack{1\leq r\leq p\\ d(s)<k\\ d(s)\leq \sigma(r)}}
X_{is}(\g) C_{sr}(\g) \tilde\l_r\,.
\end{eqnarray}
The first term of (\ref{2split}) can be written as 
\begin{eqnarray*}
\sum_{\substack{1\leq r\leq p\\ k\leq \sigma(r)}} C_{ir}(\g) \tilde\l_r &=&
\sum_{\substack{1\leq r\leq p\\ k=\sigma(r)}} C_{ir}(\g) \tilde\l_r +
\sum_{\substack{1\leq r\leq p\\ k< \sigma(r)}}
C_{ir}(\gamma) \tilde\l_r.
\end{eqnarray*}
From \eqref{matriceC}, the Lipschitz function $C_{ir}(x)$
vanishes at zero when $\alpha_k+1\leq j\leq m_k$ and $d(i)=\sigma(r)$,
then $C_{ir}(\gamma(t))=O(t)$ and 
\begin{eqnarray}\label{stima2.1}
\sum_{\substack{1\leq r\leq p\\ k\leq \sigma(r)}} C_{ir}(\g) \tilde\l_r
=\sum_{\substack{1\leq r\leq p\\ k= \sigma(r)}} O(t)\; t^{k-1}
+ \sum_{\substack{1\leq r\leq p\\ k< \sigma(r)}} O(1)\;t^{\sigma(r)-1}
= O(t^k).
\end{eqnarray}
Let us now consider the second term of \eqref{2split}.
According to (\ref{hompolprof}), we know that $X_{is}(\gamma(t))=O(t^{k-d(s)})$.
Unfortunately, this estimate is not enough for our purposes,
as one can check observing that $\tilde\lambda_r=O(t^{\sigma(r)-1})$
and $C_{sr}=O(1)$ for some of $s,r$.
To improve the estimate on $X_{is}$ we will use Lemma~\ref{subalgebra},
according to which the subspace spanned by 
$$
\left(X_1^1,\ldots,X^1_{\alpha_1},\ldots,X^\iota_1,\ldots,
X_{\alpha_\iota}^\iota\right)
$$
is a subalgebra. Then we define
$$
\mathcal F =\mbox{\rm span}\{X^k_s\mid 1\leq k\leq\iota\,,\,
1\leq s\leq\alpha_k\}
$$
along with the set $J$, that is given by the condition
$$
\mathcal F=\mbox{\rm span}\{X_j:j\in J\}.
$$
We first notice that $i\notin J$, due to our assumption
$\alpha_k+1\leq j\leq m_k$. This will allow us to apply 
Lemma~\ref{campisottoalg}, according to which we have
$$
P_i(x,y)=x_i + y_i + Q_i(x,y)=x_i + y_i + \sum_{l\notin J,\:d(l)<k} (x_l R_{il}(x,y) + y_l S_{il}(x,y))\,.
$$
As a result, assuming that $s\in J$, we obtain the key formula
$$
X_{is}(x)=\frac{\partial P_i}{\partial y_s}(x,0) = \sum_{l\notin J,\:d(l)<k} x_l\; \frac{\partial R_{il}}{\partial y_s}(x,0)\,,
$$
where $\partial_{y_s}R_{il}(x,0)$ is a homogeneous polynomial
of degree $k-d(s)-d(l)$.
By both inductive hypothesis and definition of $J$, we get
$$
\gamma_l(t)=O(t^{d(l)+1})\,,
$$ 
for every $l\notin J$ such that $d(l)<k$.
By these estimates, we achieve
\begin{multline*}
X_{is}(\gamma(t))=\sum_{l\notin J,\:d(l)<k} \gamma_l(t)
\,\frac{\partial R_{il}}{\partial y_s}(\g(t),0)
=\sum_{l\notin J,\:d(l)<k} O(t^{d(l)+1}) O(t^{k-d(s)-d(l)})
= O(t^{k+1-d(s)})\,.
\end{multline*}
Then it is convenient to split the second term of \eqref{2split}
as follows
\begin{multline}\label{stima2}
\sum_{\substack{r=1,\dots,p\\ d(s)< k\\ d(s)\leq \sigma(r)}}
X_{is}(\gamma)\, C_{sr}(\g) \tilde\l_r
 = \sum_{\substack{r=1,\dots,p\\ d(s)< k\\ d(s)\leq \sigma(r)\\ s\in J}}
X_{is}(\gamma)\, C_{sr}(\g)\, \tilde\lambda_r +
\sum_{\substack{r=1,\dots,p\\ d(s)< k\\ d(s)\leq \sigma(r)\\ s\notin J}} X_{is}(\gamma)\, C_{sr}(\g)\, \tilde\lambda_r\,,
\end{multline}
where the first addend of the previous decomposition
can be estimated as
\begin{eqnarray}
\sum_{\substack{1\leq r\leq p\\ d(s)< k\\ d(s)\leq \sigma(r)\\ s\in J}}
X_{is}(\gamma)\, C_{sr}(\gamma)\, \tilde\l_r &=&
\sum_{\substack{1\leq r\leq p\\ d(s)< k\\ d(s)\leq \sigma(r)\\ s\in J}} 
O(t^{k+1-d(s)})\,O(1)\,O(t^{\sigma(r)-1})\,=\, O(t^k)\,.\label{stima2.2}
\end{eqnarray}
Finally, we consider the second addend of \eqref{stima2}, writing
it as the following sum
\begin{equation}\label{lasterm}
\sum_{\substack{1\leq r\leq p\\ d(s)< k\\ d(s)\leq \sigma(r)\\ s\notin J}} X_{is}(\gamma)\, C_{sr}(\gamma)\, \tilde\l_r
= \sum_{\substack{1\leq r\leq p\\ d(s)< k\\ d(s)= \sigma(r)\\ s\notin J}} 
X_{is}(\gamma)\,C_{sr}(\gamma)\,\tilde\lambda_r +
\sum_{\substack{1\leq r\leq p\\ d(s)< k\\ d(s)< \sigma(r)\\ s\notin J}} 
X_{is}(\gamma)\, C_{sr}(\gamma)\, \tilde\lambda_r\,.
\end{equation}
The first term of \eqref{lasterm} can be written as
\begin{eqnarray*}
\sum_{\substack{1\leq r\leq p\\ d(s)< k\\ d(s)=\sigma(r)\\ s\notin J}}
O(t^{k-d(s)})\,O(t)\, O(t^{\sigma(r)-1})=O(t^k)\,,
\end{eqnarray*}
where we have used the fact that
$C_{sr}(x)=O(|x|)$ when $d(s)= \sigma(r)$ and $s\notin J$,
according to \eqref{matriceC}. The second term of \eqref{lasterm}
corresponds to the sum
\begin{eqnarray*}
 \sum_{\substack{1\leq r\leq p\\ d(s)< k\\ d(s)< \sigma(r)\\ s\notin J}} 
O(t^{k-d(s)})\, O(1)\, O(t^{\sigma(r)-1})
= O(t^k)\,.
\end{eqnarray*}
As a result, the second term of \eqref{stima2} is also equal
to some $O(t^k)$, hence taking into account \eqref{stima2.2}
we get that the second term of \eqref{2split} is $O(t^k)$.
Thus, taking into account \eqref{2split} and \eqref{stima2.1}
we achieve $\dot\gamma(t)=O(t^k)$, which proves the second part of \eqref{stimesopra} and ends the proof.
\end{proof}
\begin{rem}\label{Ounif}
{\upshape Analyzing the previous proof, it is easy to realize that
the functions $O(t^k)$ appearing in the statement of Lemma~\ref{curve}
can be estimated by $t^k$, uniformly with respect to $\lambda$ varying
in a compact set: there exists a constant $M>0$ such that
\begin{equation}\label{uniformita}
\begin{array}{ll}
\left| \gamma_j^k(t,\lambda) - \bigl[\l_j^k/k + g_j^k(\l_1^1,\dots,\lambda_{\alpha_{k-1}}^{k-1})
\bigr] t^k \right| \leq Mt^{k+1} & \text{if }1\leq j\leq \alpha_k\vspace{.1cm}\\
|\gamma_j^k(t,\lambda)| \leq Mt^{k+1} & \text{if }\alpha_k+1\leq j\leq m_k.
\end{array}
\end{equation}
for all $\lambda$ belonging to a compact set $L$ and every $t<t_0$: here and in the following, we have set $\g_\l:=\g(\cdot,\l)$.}
\end{rem}
Our next step will be to prove that our curves $\gamma(\cdot,\lambda)$
cover a neighbourhood of a point with maximum degree.
To do this, we fix graded coordinates with respect to the basis
$(X^k_j)$ and consider the diffeomorphism
$G:\R^p\longrightarrow\R^p$
arising from Lemma~\ref{curve} and that can be associated
with any point of maximum degree in a $C^{1,1}$ smooth submanifold.
We set
\begin{eqnarray}\label{defg}
G_i(\lambda)=\lambda_i/\sigma(i)+g_i(\lambda_1,\ldots,
\lambda_{\sum_{s=1}^{\sigma(i)-1}\alpha_s})\,,
\end{eqnarray}
where $(g_1,\ldots,g_p)
=(g_1^1,\dots,g^1_{\alpha_1},\dots,g^\iota_1,\dots,
g_{\alpha_\iota}^\iota)$
and $g^k_j$ are given by Lemma~\ref{curve}.
Then $G(0)=0$ and by explicit computation of the inverse function,
the definition (\ref{defg}) implies global invertibility of $G$.
\begin{rem}\label{remG}{\rm
The diffeomorphism $G$ also permits us to state Lemma~\ref{curve}
as follows
\begin{eqnarray}\label{Pdifc}
	\gamma(t,\lambda)=\delta_t\left(G(\lambda)+O(t)\right)\in\R^q\,,
\end{eqnarray}
where $G(\lambda)$ belongs to $\R^p\times\{0\}$, precisely,
it lies in the $p$-dimensional subspace
$\Pi_\Sigma(x)$ with respect to the associated graded coordinates.
}\end{rem}
We will denote by $c(t,\lambda)$ the projection of
$\gamma(t,\lambda)$ on $\Pi_\Sigma(x)$, namely
\begin{eqnarray}\label{clambda}
c(t,\lambda)=\tilde\pi_\Sigma(x)\big(\gamma(t,\lambda)\big),
\end{eqnarray}
where $\tilde\pi_\Sigma(x)$ represents $\pi_\Sigma(x)$
of \eqref{canproj} with respect to graded coordinates
arising from \eqref{basispi}.
In the sequel, the estimates
\begin{eqnarray}\label{ci}
c_i(t,\lambda)=G_i(\l)t^{\sigma(i)} + O(t^{\sigma(i)+1})
\end{eqnarray}
will be used. They follow from Lemma~\ref{curve} and the
definitions of $c$ and $G$.
\begin{lem}\label{lemt0}
There exists $t_0>0$ such that for every
$t_1\in]0,t_0[$, there exists a neighbourhood $V$
of $0$ such that
$$
V\cap\Sigma\subset\left\{ \gamma(t,\lambda): \lambda\in G^{-1}(S^{p-1})
\text{ and }0\leq t<t_1 \right\}.
$$
\end{lem}
\begin{proof}
We fix $t_0>0$ as in Lemma~\ref{curve}, where we have chosen
$L=G^{-1}(S^{p-1})$. Let $t_1\in]0,t_0[$ be arbitrarily fixed.
Taking into account
Corollary~\ref{graficosulpiano}, it suffices to prove that
the set $\{ c(t,\lambda): \l\in L,0\leq t<t_1\}$ covers a neighbourhood
of 0 in $\R^p$. 
For each $t\in]0,t_1[$, we define the ``projected dilations''
$\Delta_t=\tilde\pi_\Sigma(x)\circ\delta_t$
corresponding
to the following diffeomorphisms of $\R^p$ 
$$
\Delta_t(y_1,\dots,y_p)=
\left(t^{\sigma(1)}y_1,\dots,t^{\sigma(i)}y_i,\dots,
t^{\sigma(p)}y_p\right).
$$
Now we can rewrite \eqref{ci} as 
\begin{eqnarray}\label{cex}
c(t,\lambda)=\Delta_t\big(G(\lambda)+O(t)\big),
\end{eqnarray}
where $O(t)$ is uniform with respect to $\lambda$
varying in $G^{-1}(S^{p-1})$, according to
Remark~\ref{Ounif}. Then we define the mapping
\begin{eqnarray*}
L_t &:& S^{p-1}\to\R^{p}\\
&& u\longmapsto \Delta_{1/t}\bigl(c\big(t,G^{-1}(u)\big)\bigr)
\end{eqnarray*}
and \eqref{cex} implies 
$$
L_t(u)=u+O(t).
$$
As a consequence, $L_t\to Id_{S^{p-1}}$ as $t\to0$, uniformly
with respect to $u$ varying in $S^{p-1}$.
Then, for any sufficiently small $0<\tau<t_1$, we have
$L_{\tau}(S^{p-1})\cap B_{1/2}=\emptyset$ and $L_{\tau}$
is homotopic to $Id_{S^{p-1}}$ in $\R^p\setminus \{A\}$ for all
$A\in B_{1/2}$.
In particular, since $Id_{S^{p-1}}$ is not homotopic to a constant,
$L_\tau$ is not homotopic to a constant in 
$\R^p\setminus \{A\}$ for all $A\in B_{1/2}$.
Now, we are in the position to prove that
$$
\left\{ c(t,\lambda): \l\in G^{-1}(S^{p-1})\text{ and }0\leq t<\tau \right\}
$$
covers the open neighbourhood of 0 in $\R^p$ given by
$\Delta_\tau(F^{-1}(B_{1/2})\cap\Pi_\Sigma(e))$
that leads us to the conclusion.
By contradiction, if this were not true, then
we could find a point $A\in B_{1/2}$ such that
$A\neq \Delta_{1/\tau}(c_\lambda(t))$ for all
$\lambda\in G^{-1}(S^{p-1})$ and $0\leq t<\tau$, but then 
\begin{eqnarray*}
H &:& [0,\tau]\times S^{p-1} \to\R^p\setminus\{A\}\\
&& (s,u)\longmapsto \Delta_{1/\tau}\bigl(c\big(s,G^{-1}(u)\big)\bigr)
\end{eqnarray*}
would provide a homotopy in $\R^p\setminus\{A\}$ between the constant 0 and $L_\tau$, which cannot exist.
\end{proof}
As important consequence of Lemma~\ref{curve},
we are in the position to give the 
\vskip.25truecm
\noindent
{\em Proof of Theorem~\ref{blowup}.}
We first notice that $\Pi_\Sigma(x)$ is a subgroup of $\G$,
due to Lemma~\ref{subalgebra}. Setting $\Sxr:=\delta_{1/r}(x^{-1}\Sigma)$, it is sufficient to prove (see \cite{AmbTil}, Proposition 4.5.5) that $\Sxr\cap D_R$ converges to
$\Pi\cap D_R$ in the Kuratowski sense, i.e. that
\begin{itemize}
\item[(i)] if $y=\lim_{n\to\infty} y_n$ for some sequence $\{y_n\}$ such that $y_n\in \Sigma_{x,r_n}\cap D_R$ and $r_n\to 0$, then
$y\in\Pi_\Sigma(x)\cap D_R$;
\item[(ii)] if $y\in\Pi_\Sigma(x)\cap D_R$, then there are
$y_r\in\Sxr\cap D_R$ such that $y_r\to y$.
\end{itemize}
It is not restrictive assuming that $x=e$.
To prove (i), we set $z_n=\delta_{r_n}(y_n)\in\Sigma\cap D_{r_nR}$.
From \eqref{Pdifc}, we can find $t_1>0$ arbitrarily small such that
\begin{equation}\label{posa}
\inf_{u\in S^{p-1}}|u+O(t_1)|>0,
\end{equation}
where $|\cdot|$ is the Euclidean norm and $O(t)$ is defined in \eqref{Pdifc}.
Then for $n$ sufficiently large and taking $t_1<t_0$
Lemma~\ref{lemt0} yields a sequence $\{\tau_n\}\subset]0,t_1[$
and $\lambda_n\in G^{-1}(S^{p-1})$ such that $\gamma(\tau_n,\lambda_n)=\delta_{r_n}y_n$. Due to \eqref{Pdifc}, we achieve
$$
\delta_{\tau_n/r_n}\left(G(\lambda_n)+O(\tau_n)\right)=y_n\,,
$$
hence \eqref{posa} implies that $\tau_n/r_n$ is bounded.
Up to subsequences, we can assume that 
$G(\lambda_n)\to\zeta$ and $\tau_n/r_n\to s$,
then $y_n\to\delta_s\zeta=y$.
From Remark~\ref{remG}, we know that $G(\lambda)\in\Pi_\Sigma(x)$
with respect to our graded coordinates, hence $y\in\Pi_\Sigma(x)$.
To prove (ii), we choose $y\in\Pi_\Sigma(x)\cap D_R$
and set $\lambda=G^{-1}(y)$. By Lemma~\ref{curve} there
exists $r_0>0$ depending on the compact set
$G^{-1}( D_R\cap\Pi_\Sigma(x))$ such that
the solution  $r\to\gamma(r,\lambda')$ of \eqref{gammalambda}
is defined on $[0,r_0]$
for every $\lambda'\in G^{-1}( D_R\cap\Pi_\Sigma(x))$.
Clearly, $\gamma(r,\lambda')\in\Sigma$, then \eqref{Pdifc}
implies that
$$
\delta_{1/r}\big(\Sigma\big)
\ni y_r=\delta_{1/r}\big(\gamma(r,\lambda)\big)
\longrightarrow G(\lambda)=y\,.
$$
This ends the proof. \hskip11.9truecm $\Box$
\vskip.25truecm

Next we prove Theorem~\ref{densita}.
We will denote by $\mbox{vol}_{\tilde g}^p$ the Riemannian $p$-dimensional
volume with respect to an arbitrary metric $\tilde g$.
\vskip.25truecm
\noindent
{\em Proof od Theorem~\ref{densita}.}
Without loss of generality we assume that
$x$ is the identity element $e$ and consider
graded coordinates $F:\R^q\longrightarrow\G$ centered at $0$ with respect to $X_j^k$. Notice that balls $F^{-1}(B_{x,r})$ in $\R^q$ through graded coordinates will be simply denoted by $B_{x,r}$. 
According to Corollary~\ref{graficosulpiano}, we parametrize
$\Sigma$ by the $C^{1,1}$ function
$\varphi:A\subset\Pi_\Sigma(e)\to\R^{q-p}$,
such that $\Sigma$ is the image of
\begin{eqnarray*}
\F &:& A\subset\Pi_\Sigma(e)\longrightarrow\R^q\\
&& y\mapsto (y_1^1,\dots,y_{\alpha_1}^1,\f_{\alpha_1+1}^1(y),\dots,\f_{m_1}^1(y),\dots,y_1^\iota,\dots,y_{\alpha_\iota}^\iota,\f_{\alpha_\iota+1}^\iota(y),\dots,\f_{m_\iota}^\iota(y)).
\end{eqnarray*}
For any sufficiently small $r>0$, we have
\begin{eqnarray}\label{densitaddim}
\lim_{r\downarrow 0}
\frac{\mbox{vol}_{\tilde g}^p(\Sigma\cap B_r)}{r^d} &=& \frac{1}{r^d}\int_{\F^{-1}(B_r)}J_{\tilde g} \F(y)\,dy \nonumber\\
&=& \int_{\Delta_{1/r}(\F^{-1}(B_r))} J_{\tilde g}\F(\Delta_r(y)) dy\,,
\end{eqnarray}
where $\Delta_r={\delta_r}_{|\Pi_{\Sigma}(e)}$ and its jacobian
is exactly equal to $d$.
Notice that $\Delta_{1/r}(\F^{-1}(B_r)) = (\delta_{1/r}\circ\F\circ\Delta_r)^{-1}(B_1)$ is the set of elements
$y\in\Pi_\Sigma(e)$ such that 
$$
\left(y_1^1,\dots,y_{\alpha_1}^1,\frac{\f_{\alpha_1+1}^1(\Delta_ry)}{r},\dots,\frac{\f_{m_1}^1(\Delta_ry)}{r},\dots,y_1^\iota,\dots,y_{\alpha_\iota}^\iota, \frac{\f_{\alpha_\iota+1}^\iota(\Delta_ry)}{r^\iota}, \dots,\frac{\f_{m_\iota}^\iota(\Delta_ry)}{r^\iota}\right)
$$
belongs to $B_1$ and that
$$
\Delta_{1/r}(\F^{-1}(B_r))=\tilde\pi_\Sigma(e)(\Sigma_{0,r}\cap B_1),
$$
where $\tilde\pi_\Sigma(e)$ is the projection $\pi_\Sigma(e)$
with respect to graded coordinates, i.e. the mapping
$$
\R^q\ni(z_1^1,\dots,z_{m_1}^1,\dots,z^\iota_1,\dots,z_{m_\iota}^\iota)\longmapsto (z_1^1,\dots,z_{\alpha_1}^1,\dots,z^\iota_1,\dots,z_{\alpha_\iota}^\iota)\in\Pi_\Sigma(e).
$$
We will denote the projection $\tilde\pi_\Sigma(e)$ by
$\pi$. By continuity of $\pi$, for every $\ep>0$ we can find a neighbourhood $\mathcal N\subset\R^q$ of $\Pi_\Sigma(e)\cap  D_1$
such that
$\pi(\mathcal N)\subset\Pi_\Sigma(e)\cap B_{1+\ep}$;
by Theorem~\ref{blowup} and the definition of Hausdorff convergence, for sufficiently small $r$ we have $\Sigma_{0,r}\cap D_1\subset\mathcal N$ and so
\begin{equation}\label{contenimento1}
\Delta_{1/r}(B_r)\subset\pi(\Sigma_{0,r}\cap D_1)
\subset \Pi_\Sigma(e)\cap B_{1+\ep}.
\end{equation}
If we also prove that 
\begin{equation}\label{contenimento2}
\Pi_\Sigma(e)\cap B_{1-\ep}\subset\Delta_{1/r}(\F^{-1}(B_r))
\end{equation}
for small $r$, we will have $\chi_{\delta_{1/r}(\F^{-1}(B_r))}\to\chi_{\Pi_\Sigma(e)\cap B_1}$ in $L^1(\Pi_\Sigma(e))$. This fact and \eqref{densitaddim} imply that
$$
\lim_{r\downarrow 0} \frac{\mbox{vol}_{\tilde g}^p(\Sigma\cap B_r)}{r^d} = J_{\tilde g}\F(0)\; \mathcal L^p (\Pi_\Sigma(e)\cap B_1) = J_{\tilde g}\F(0)\, \theta(\tau_{\Sigma}^d(0)).
$$
By Corollary~\ref{graficosulpiano} we know that $\nabla\F(0)=C(0)$, where $C$ is given by Lemma \ref{buonivettoritangenti}; therefore $J_{\tilde g}\F(0)$ must coincide with the Jacobian of the matrix $C(0)$, i.e. with $|v_1(0)\wedge\dots\wedge v_p(0)|_{\tilde g}$.
By virtue of Corollary~\ref{vettangentemax}, we have 
$$
|{\tau_{\Sigma}^d}(e)|=\left|\frac{X_1^1\wedge\dots\wedge X_{\alpha_1}^\iota\wedge\dots\wedge X^\iota_1\wedge\dots\wedge X_{\alpha_\iota}^\iota}
{|v_1(0)\wedge\dots \wedge v_p(0)|_{\tilde g}}\right|_g
=\frac{1}{|v_1(0)\wedge\dots\wedge v_p(0)|_{\tilde g}}.
$$
Finally, it remains to prove \eqref{contenimento2}.
We fix
$$
y=(y_1,\dots,y_p)=
(y_1^1,\dots,y^1_{\alpha_1},\dots,y^\iota_{\alpha_\iota})
\in\Pi_\Sigma(e)\cap B_{1-\ep}
$$
and set $z:=\delta_r(y)\in B_{(1-\ep)r}$.
Let $t_0>0$ be as in Lemma~\ref{lemt0} and consider
$t_1\in]0,t_0[$ to be chosen later.
By the same lemma, for every $r>0$ sufficiently small
there exist $\lambda\in G^{-1}(S^{p-1})$ and $t\in[0,t_1[$
such that $\F(z)=\gamma(t,\lambda)$.
Since $|G(\l)|=1$, we can find $1\leq i\leq p$
such that $|G_i(\lambda)|\geq 1/\sqrt{p}$.
Notice that
\begin{eqnarray}\label{pic}
\pi_\Sigma(e)(\Phi(z))=z=\pi_\Sigma(e)(\gamma(t,\lambda))=c(t,\lambda),
\end{eqnarray}
then \eqref{ci} implies
$$
Mt^{\sigma(i)+1}\geq|G_i(\lambda)|t^{\sigma(i)}-|z_i|
\geq t^{\sigma(i)}/\sqrt{p}-|y_i|r^{\sigma(i)},
$$
where $M>0$ is given in Remark~\ref{Ounif}
with $L=G^{-1}(S^{p-1})$. It follows that
$$
(1/\sqrt{p}-Mt_1)t^{\sigma(i)}\leq (1/\sqrt{p}-Mt)t^{\sigma(i)}
\leq |y_i|\,r^{\sigma(i)}\,.
$$
Now, we can choose $t_1>0$ such that $1/\sqrt{p}-Mt_1\geq \epsilon>0$,
getting a constant $N>0$ depending only on $p$, $|y|$ and $M$ such that
\begin{eqnarray}\label{constN}
t\leq N\,r\,.
\end{eqnarray}
Taking into account \eqref{pic}
and the explicit estimates of \eqref{uniformita}, we get
some $1\leq k\leq\iota$ and $\alpha_j+1\leq j\leq m_j$ such that
$$
|c_i(t,\lambda)|=|\gamma^k_j(t,z)|=|\varphi^k_j(z)|
\leq M t^{k+1}\,,
$$
where we notice that $k=\sigma(i)$. By \eqref{constN}, the previous
estimate yield
\begin{eqnarray}\label{Mrunif}
|\varphi^k_j(\delta_ry)|=|\varphi^k_j(z)|\leq \tilde M r^{k+1}\,,
\end{eqnarray}
where $\tilde M=M N^{k+1}$. Estimate \eqref{Mrunif} has been
obtained with $\tilde M$ independent from $r>0$ sufficiently small.
Therefore 
\begin{eqnarray*}
\left(y_1^1,\dots,y_{\alpha_1}^1,\frac{\f_{\alpha_1+1}^1(\delta_ry)}{r},\dots,\frac{\f_{m_1}^1(\delta_ry)}{r},\dots,
y_1^\iota,\dots,y_{\alpha_\iota}^\iota, \frac{\f_{\alpha_\iota+1}^\iota(\delta_ry)}{r^\iota}, \dots,\frac{\f_{m_\iota}^\iota(\delta_ry)}{r^\iota}\right)
\end{eqnarray*}
belongs to $B_1$ definitely as $r$ goes to zero, namely,
$y\in \Delta_{1/r}\Phi^{-1}(B_r)$ for $r>0$ small enough.
We observe that $N$ linearly depends on $|y|$ and is independent
from $r>0$, then the constant $\tilde M$ in \eqref{Mrunif}
can be fixed independently from $y$ varying in the
bounded set $\Pi_\Sigma(e)\cap B_{1-\ep}$,
whence \eqref{contenimento2} follows.
\hskip14truecm $\Box$
\begin{rem}\label{coin1cod}{\rm
In the case $x$ is a non-horizontal point of a $p$-dimensional submanifold $\Sigma$, namely $d_\Sigma(x)=Q$$-$$k$,
the limit \eqref{eqdensita} coincides with the limit
stated in Theorem~3.5 of \cite{Mag8}, according to which
\begin{equation}\label{blwlim1}
\frac{\mbox{vol}_{\tilde{g}}(\Sigma\cap B_{p,r})}{r^{Q-k}}
\longrightarrow
c(\tilde g,g)\;
\frac{\theta\big({\bf{\tilde n}}_{g,H}(p)\big)}
{|{\bf{\tilde n}}_{g,H}(p)|}\,.
\end{equation}
Here we have defined 
$$
c(\tilde g,g)=\frac{\mu_{\tilde g}(B_1)}{\mu_g(B_1)}
=\frac{\det\big(\tilde g_{ij}\big)}{\det\big(g_{ij}\big)}
$$
where $\mu_{\tilde g}$ and $\mu_g$ are the Riemannian volume measures on
$\G$ with respect to $\tilde g$ and $g$. In addition,
we have assumed that $\mu_{\tilde g}$ is also an Haar measure of $\G$. 

To show this fact, we fix a basis
$\{X_\alpha\}_{\alpha\in I}$ of
$\Lambda_k(V_1)$, where
$$X_\alpha=X_{\alpha_1}\wedge\cdots\wedge X_{\alpha_k}\quad\mbox{and}\quad
I=\{(\alpha_1,\ldots,\alpha_k)\mid
1\leq \alpha_1<\cdots<\alpha_k\leq m\}$$
and $m=\dim V_1$.
Notice that existence of a non-horizontal point
implies the condition $m\geq k$.
We choose an orthonormal frame $(B_1,\ldots,B_q)$
with respect to $\tilde g$ and define $B=B_1\wedge\cdots B_q$.
By definition of Hodge operator with respect to $\tilde g$ and
orientation $B$, we have
\begin{eqnarray}\label{keyhodge}
X_\alpha\wedge\tau_\Sigma(x)
=(-1)^{kp}\;\left\lan X_\alpha,
\tilde *\big(\tau_\Sigma(x)\big)\right\ran_{\tilde g }\; B
=(-1)^{kp}\;\left\lan X_\alpha,
{\bf \tilde n}(x)\right\ran_{\tilde g }\; B\,.
\end{eqnarray}
Taking into account the decomposition
$$
\tau_\Sigma(x)=\sum_{\alpha\in I}c_{*\alpha}\; *\big(\!X_\alpha\big)+
\sum_{d(\beta)<Q-k} c_\beta X_\beta\,,
$$
where $\alpha\cup*\alpha=\{1,\ldots,q\}$, we have
$$
\tau^d_\Sigma(x)=\sum_{\alpha\in I}c_{*\alpha} *X_\alpha
$$
where $d=Q-k$. In addition, we have
$B=c(g,\tilde g)\,X_1\wedge\cdots\wedge X_q$.
Then \eqref{keyhodge} implies
$$
c_{*\alpha}\,X_\alpha\wedge X_{*\alpha}
=c(g,\tilde g)\,(-1)^{kp}\;\left\lan X_\alpha,
{\bf \tilde n}(x)\right\ran_{\tilde g }\; X_1\wedge\cdots\wedge X_q\,.
$$
According to Definition~2.6 of \cite{Mag8}, we have
$
\left\lan X_\alpha,{\bf \tilde n}(x)\right\ran_{\tilde g }
=\left\lan X_\alpha,{\bf \tilde n}_g(x)\right\ran
$
and
\begin{eqnarray}
{\bf \tilde n}_{g,H}(x)
=\sum_{\alpha\in I}
\left\lan X_\alpha,{\bf \tilde n}_g(x)\right\ran X_\alpha\,,
\end{eqnarray}
therefore
$$
|\tau^d_\Sigma(x)|^2=\sum_{\alpha\in I}(c_{*\alpha})^2
=c(g,\tilde g)^2\,|{\bf \tilde n}_{g,H}(x)|^2,
$$
concluding our proof.
}\end{rem}
As it has been mentioned in the introduction,
it is easy to find groups where non-horizontal submanifolds
of a given topological dimension cannot exist.
\begin{ex}\label{exanothor}{\rm
Let us consider the 5-dimensional stratified group
$\mathbb E^5$ with a basis $X_1,\ldots, X_5$ subject to the
only nontrivial relations 
$$
[X_1,X_2]=X_3,\quad[X_1,X_3]=X_4,\quad[X_1,X_4]=X_5
$$
and the grading
$$
V_1=\mbox{span}\{X_1,X_2\},\quad V_2=\mbox{span}\{X_3\},\quad
V_3=\mbox{span}\{X_4\},\quad V_4=\mbox{span}\{X_5\}\,.
$$
Then $m=2$ and a 2-dimensional submanifold has codimension $k=3$.
As a result, $m-k<0$ hence any 2-dimensional submanifold $\Sigma$
satisfies $d(\Sigma)<Q-k=11-3=8$. In other words, all 2-dimensional
submanifolds of $\mathbb E^5$ are horizontal.
}\end{ex}

\section{Some applications in the Engel group}\label{engelsbman}

In this section we wish to present examples of
2-dimensional submanifolds of all possible degrees
in the Engel group $\mathbb E^4$.

We represent $\mathbb E^4$ as $\mathbb R^4$ equipped
with the vector fields
$
X_i=\sum_{j=1}^4 A_i^j(x)e_j\,,
$
where 
$$
A(x)=\left(\begin{array}{cccc}
1 & 0              &     0   & 0 \\
0 & 1              &     0   & 0 \\
0 & x_1            &     1   & 0 \\
0 & x_1^2/2        &     x_1 & 1
\end{array}\right)\,,
$$
$(e_1,e_2,e_3,e_4)$ is the canonical basis of $\mathbb R^4$
and $x=(x_1,x_2,x_3,x_4)$.

Let $\Phi:U\longrightarrow{\mathbb R}^4$ be the parametrization
of a 2-dimensional submanifold $\Sigma$,
where $U$ is an open subset of $\mathbb R^2$.
We set $u=(u_1,u_2)=(x,y)\in U$ and consider 
$\Phi_{u_i}=\sum_{j=1}^4\Phi_{u_i}^j\,e_j$. Taking into account that
$$
A(x)^{-1}=\left(\begin{array}{cccc}
1 & 0              &     0   & 0 \\
0 & 1              &     0   & 0 \\
0 & -x_1           &     1   & 0 \\
0 & x_1^2/2        &   - x_1 & 1
\end{array}\right)
$$
and that
\begin{eqnarray}\label{ei}
e_i=\sum_{j=1}^4(A(x)^{-1})_i^j X_j\,,
\end{eqnarray}
we obtain
$$
\Phi_{u_i}=\Phi_{u_i}^1 X_1+\Phi_{u_i}^2 X_2+
\left(\Phi_{u_i}^3-\Phi^1\,\Phi_{u_i}^2\right) X_3+
\left(\Phi_{u_i}^4-\Phi^1\,\Phi_{u_i}^3+\frac{(\Phi^1)^2}{2}\,
\Phi_{u_i}^2\right) X_4\,.
$$
It follows that
\begin{eqnarray}
&&\Phi_x\wedge\Phi_y=
\Phi_u^{12} X_1\wedge X_2+\left(\Phi_u^{13}-\Phi^1\Phi_u^{12}\right)
X_1\wedge X_3
+\left(\Phi_u^{14}-\Phi^1\,\Phi_u^{13}+\frac{(\Phi^1)^2}{2}\,
\Phi_u^{12}\right)X_1\wedge X_4
\label{phixy} \\
&&+\Phi_u^{23}  X_2\wedge X_3
+\left(\Phi_u^{24}-\Phi^1\,\Phi_u^{23}\right) X_2\wedge X_4
+\left(\Phi_u^{34}+\frac{(\Phi^1)^2}{2}\Phi_u^{23}-\Phi^1\Phi_u^{24}\right) X_3\wedge X_4\,,
\nonumber
\end{eqnarray}
where we have set
$$
\Phi_u^{ij}=\det\left(\begin{array}{cc}
\Phi_x^i & \Phi_y^i \\
\Phi_x^j & \Phi_y^j
\end{array}\right)\,.
$$
In the sequel, we will use \eqref{phixy} to obtain nontrivial
examples of 2-dimensional submanifolds with different degrees
in $\mathbb E^4$.

\begin{rem}{\rm
Recall that 2-dimensional submanifolds of degree 2 in $\mathbb E^4$
cannot exist, due to non-integrability of the horizontal distribution
$\mbox{span}\{X_1,X_2\}$.
}\end{rem}
The next example wants to give a rather general method to obtain
nontrivial examples of 2-dimensional submanifolds of degree 3. 
Clearly, the submanifold
$\{(0,x_2,x_3,0\}$ is the simplest example, as one can check
using \eqref{phixy}.
\begin{ex}{\rm
Having degree three means that the first order
fully non-linear conditions
\begin{eqnarray}\label{flsyst}
\left\{\begin{array}{l}
\Phi_u^{34}+\frac{(\Phi^1)^2}{2}\Phi_u^{23}-\Phi^1\Phi_u^{24}=0 \\
\Phi_u^{24}-\Phi^1\,\Phi_u^{23}=0 \\
\Phi_u^{14}-\Phi^1\,\Phi_u^{13}+\frac{(\Phi^1)^2}{2}\,\Phi_u^{12}=0
\end{array}\right.
\end{eqnarray}
must hold.
By elementary properties of determinants, one can realize that
the previous system is equivalent to requiring that
\begin{eqnarray}
&&\nabla\Phi^3\quad\quad\mbox{is parallel to}\quad\quad
\nabla\Phi^4+\frac{(\Phi^1)^2}{2}\,\nabla\Phi^2\,, \label{first}\\
&&\nabla\Phi^2\quad\quad\mbox{is parallel to}\quad\quad
\nabla\Phi^4-\Phi^1\,\nabla\Phi^3\,, \label{second}\\
&&\nabla\Phi^1\quad\quad\mbox{is parallel to}\quad\quad
\nabla\Phi^4-\Phi^1\,\nabla\Phi^3+\frac{(\Phi^1)^2}{2}\nabla
\Phi^2\,. \label{third}
\end{eqnarray}
We restrict our search to submanifolds with $\Phi^1(x,y)=x$
and $\Phi_u^{23}\neq0$ on $U$.
Then \eqref{first} and \eqref{second}
yield functions $\lambda,\mu:U\longrightarrow\R$ such that
\begin{eqnarray*}
&&\nabla\Phi^4+\frac{x^2}{2}\,\nabla\Phi^2=\lambda\;\nabla\Phi^3 \\
&&\nabla\Phi^4-x\,\nabla\Phi^3=\mu\;\nabla\Phi^2\,.
\end{eqnarray*}
Due to the previous equations and the condition $\Phi_u^{23}\neq0$,
we impose the assumptions
$$
\lambda(u)=x\qquad\mbox{and}\qquad\mu(u)=-\frac{x^2}{2}\,.
$$
It follows that
\begin{eqnarray}\label{fineq}
\nabla\Phi^4=-\frac{x^2}{2}\nabla\Phi^2+x\,\nabla\Phi^3\,,
\end{eqnarray}
then all conditions \eqref{first}, \eqref{second} and
\eqref{third} are satisfied, namely, the system \eqref{flsyst} holds
whenever we are able to find $\Phi^4$ satisfying \eqref{fineq}.
Notice that condition \eqref{fineq} characterizes all
2-dimensional submanifolds of degree three under the condition
$\Phi_u^{23}\neq0$.
Clearly, we have an ample choice of families of functions
$\Phi^2,\Phi^3,\Phi^4$ satisfying \eqref{fineq}.
We choose the injective embedding 
of $\mathbb R^2$ into $\mathbb R^4$ defined by
\begin{eqnarray*}
\Phi(x,y)=\left(\begin{array}{c}
x
 \\
x+e^y     \\
x e^y+\frac{x^2}{2} \\
\frac{x^3}{6}+\frac{x^2}{2}e^y 
\end{array}\right)\,.
\end{eqnarray*}
One can check that $d_\Sigma(\Phi(x,y))=3$ for every $(x,y)\in\mathbb R^2$,
where $\Sigma=\Phi(\mathbb R^2)$. Here the part of $\tau_\Sigma$
with maximum degree is 
$$
\tau_\Sigma^3\big(\Phi(x,y)\big)=
-\frac{e^{2y}}{\sqrt{e^{4y}+x^2e^{4y}+\frac{x^4}{4}e^{4y}}}
X_2\wedge X_3
$$
and due to (\ref{integrfor}), the spherical Hausdorff measure of bounded portions of $\Sigma$ is positive and finite.
}\end{ex}
It is clear that submanifolds of higher degree are easier to 
be contructed.
\begin{ex}\label{deg4}{\rm
Let us consider
$$
\Phi(x,y)=\left(x,y,\frac{y^2}{2},\frac{y^2}{2}\right).
$$ 
Then we have
\begin{eqnarray*}
&&\Phi_u^{12}=1,\qquad \Phi_u^{13}=y,\qquad \Phi_u^{14}=y  \\
&&\Phi_u^{23}=0,\qquad \Phi_u^{24}=0,\qquad \Phi_u^{34}=0.
\end{eqnarray*}
By \eqref{phixy} we have
\begin{eqnarray}
\Phi_x\wedge\Phi_y=X_1\wedge X_2+\left(y-x\right)X_1\wedge X_3
+\left(y-xy+\frac{x^2}{2}\right) X_1\wedge X_4.
\end{eqnarray}
Recall that $\Sigma_r$ is the subset of points in $\Sigma$
with degree equal to $r$. With this notation we have
\begin{eqnarray*}
&&\Sigma_4=\left\{\Phi(x,y)\;\Big|\; y\in]0,2[ \right\}\bigcup
\left\{\Phi(x,y)\;\Big|\; y\in{\mathbb R}\setminus [0,2]
\quad\mbox{and} \quad |y-x|^2\neq y^2-2y \right\} \\
&&
\Sigma_3=\left\{\Phi(y+\sigma\sqrt{y^2-2y},y)\;\Big|\;
\sigma\in\{1,-1\}\quad\mbox{and}\quad
y\in{\mathbb R}\setminus [0,2]
\right\} \\
&&\Sigma_2=\left\{\Phi(0,0),\Phi(2,2)\right\}\,.
\end{eqnarray*}
We will check that the curve
$$
{\mathbb R}\setminus [0,2]\ni y\longrightarrow \gamma(y)
=\Phi(y+\sigma\sqrt{y^2-2y},y)
$$
with $\sigma\in\{1,-1\}$ have degree constantly equal to 2.
Due to \eqref{ei}, we achieve 
$$
\dot\gamma=\dot\gamma^1 X_1+\dot\gamma^2 X_2+
\left(\dot\gamma^3-\gamma^1\,\dot\gamma^2\right) X_3+
\left(\dot\gamma^4-\gamma^1\,\dot\gamma^3+
\frac{(\gamma^1)^2}{2}\,\dot\gamma^2\right) X_4\,,
$$
where one can check that 
\begin{eqnarray}
\left(\dot\gamma^4-\gamma^1\,\dot\gamma^3+
\frac{(\gamma^1)^2}{2}\,\dot\gamma^2\right)=0\quad\mbox{and}\quad
\left(\dot\gamma^3-\gamma^1\,\dot\gamma^2\right)=-\sigma\sqrt{y^2-2y}\neq0\,.
\end{eqnarray}
It follows that $\Sigma_3$ is the union of two curves with degree
constantly equal to 2. Applying \eqref{integrfor} we get that
$\cS^2\res\Sigma_3$ is positive and finite on bounded open
pieces of $\Sigma_3$, hence $\cS^4(\Sigma_3)=0$.
In particular, we have proved that
$$
\cS^4(\Sigma\setminus\Sigma_4)=0,
$$
then the Hausdorff dimension of $\Sigma$ is 4 and furthermore
$\cS^4\res\Sigma$ is positive and finite on open bounded pieces
of $\Sigma$. Clearly, \eqref{integrfor} holds.
}\end{ex}
\begin{ex}{\rm
Using \eqref{phixy} one can check that
2-dimensional submanifolds given by
$$
\Phi(x,y)=\left(\begin{array}{c}
0 \\
\Phi^2(x,y) \\
\Phi^3(x,y) \\
\Phi^4(x,y) 
\end{array}\right)
$$
where $\Phi_u^{34}\neq0$ have degree 5$=Q-k$. Notice that
these submanifolds are then non-horizontal.
}\end{ex}
\begin{rem}\label{nrpsubgroup}{\rm
Let us consider $\Sigma$ as in Example~\ref{deg4}.
It is easy to check that
$$
\delta_{1/r}\Sigma\cap D_R\longrightarrow S\cap D_R
$$ 

where
$$
S=\{(x_1,0,0,x_4)\mid x_4\geq0\}\,.
$$
Clearly, $S$ cannot be a subgroup of $\mathbb E^4$, since
all $p$-dimensional subgroups of stratified groups
are homeomorphic to $\mathbb R^p$, see \cite{Vara}.
This fact, may occur since the origin in $\Sigma$
has not maximum degree, as one can check in
Example~\ref{deg4}.
}\end{rem}

{\bf Acknowledgements.}
We are grateful to Giuseppe Della Sala for fruitful discussions.

\end{document}